\documentclass[a4paper]{scrartcl}
\usepackage[utf8]{inputenc}
\usepackage{amsmath,amssymb,amsthm,amstext}
\usepackage[german,english]{babel}
\usepackage[alignlabels,midshaft,nohug]{diagrams}
\usepackage{hyperref}
\usepackage{multicol}
\usepackage{graphicx}

\setkomafont{sectioning}{\rmfamily\bfseries}

\newcommand{\cites}{\cite}

\newenvironment{proof2}{\begin{proof}[\textsc{Proof}]}{\end{proof}}

\newcommand{\IR}{\mathbb R}

\newcommand{\CR}{\mathcal R}

\newcommand{\rank}{\operatorname{rank}}
\newcommand{\CL}{\mathcal L}
\newcommand{\CS}{\mathcal S}

\newcommand{\IZ}{\mathbb Z}
\newcommand{\IT}{\mathbb T}
\newcommand{\IC}{\mathbb C}

\newcommand{\CC}{\mathcal C}

\newcommand{\IN}{\mathbb N}

\renewcommand{\phi}{\varphi}
\newcommand{\eps}{\varepsilon}
\newcommand{\gdw}{\Leftrightarrow}

\newcommand{\fa}{\forall}
\newcommand{\ex}{\exists\, }

\newcommand{\diag}{\operatorname{diag}}
\newcommand{\dist}{\operatorname{dist}}

\newcommand{\id}{\operatorname{id}}

\newcommand{\CM}{\mathcal M}

\newcommand{\CP}{\mathcal P}

\newcommand{\spr}[1]{\left\langle #1 \right\rangle}

\newcommand{\CF}{\mathcal F}
\newcommand{\dif}{\mathrm d}
\newcommand{\Dif}{\mathrm D}

\newcommand{\into}{\hookrightarrow}

\newcounter{punkt}
\theoremstyle{plain}
\newtheorem{lemma}[punkt]{Lemma}

\newtheorem{thm}[punkt]{Theorem}
\newtheorem{prop}[punkt]{Proposition}
\theoremstyle{remark}

\theoremstyle{definition}
\newtheorem{definition}[punkt]{Definition}

\newtheorem{assumption}[punkt]{Hypothesis}
\numberwithin{equation}{section}
\numberwithin{punkt}{section}

\title{Equilibria for the $N$--vortex--problem in a general bounded domain}
\author{Christian Kuhl}
\date{\today}

\begin{document}

  \maketitle

\begin{abstract}
This article is concerned with the study of existence and properties of stationary solutions for the dynamics of $N$ point vortices
in an idealised fluid constrained to a bounded two--dimen\-sional domain $\Omega$, which is governed by a Hamiltonian system
\[
\left\{
\begin{aligned}
\Gamma_i\frac{\dif x_i}{\dif t}
 &=\frac{\partial H_\Omega}{\partial y_i}(z_1,\dots,z_N)\\ 
\Gamma_i\frac{\dif y_i}{\dif t}
 &=-\frac{\partial H_\Omega}{\partial x_i}(z_1,\dots,z_N)
\end{aligned} 
\hspace{2cm}\text{where}\ z_i=(x_i,y_i),\ i=1,\dots,N,
\right.
\]
where $H_\Omega(z):=\sum_{j=1}^N\Gamma_j^2h(z_j)+\sum_{i,j=1, i\not=j}^N\Gamma_i\Gamma_jG(z_i,z_j)$ is the so--called Kirchhoff--Routh--path function under various conditions on the ``vorticities'' $\Gamma_i$ and various topological and geometrical assumptions on $\Omega$. 

In particular, we will prove that (under an additional technical assumption) if it is possible to align the vortices along a line, such that the signs of the $\Gamma_i$ are alternating and $|\Gamma_i|$ is increasing, $H_\Omega$ has a critical point.

If $\Omega$ is not simply connected, we are able to derive a critical point of $H_\Omega$, if $\sum_{j\in J}\Gamma_j^2>\sum_{\substack{i,j\in J\\ i\not=j}}|\Gamma_i\Gamma_j|$ for all $J\subset\{1,\dots,N\}$, $|J|\ge 2$.

\end{abstract}

  \tableofcontents

  \pagestyle{plain}

\section{Introduction}
The $N$--vortex--problem of fluid dynamics is concerned with the dynamics of
$N$ point vortices $z_1,\dots,z_N$ in an ideal fluid constrained to a two--dimensional
domain $\Omega$ with corresponding vortex strengths (so-called vorticities) $\Gamma_1,\dots,\Gamma_N\in\IR$,
whose absolute values determine the degree to which the surrounding fluid is curled and whose
signs determine the direction of revolution for the surrounding fluid. It is governed by a Hamiltonian
system
\begin{equation}\label{1}
\left\{
\begin{aligned}
\Gamma_i\frac{\dif x_i}{\dif t}
 &=\frac{\partial H_\Omega}{\partial y_i}(z_1,\dots,z_N)\\ 
\Gamma_i\frac{\dif y_i}{\dif t}
 &=-\frac{\partial H_\Omega}{\partial x_i}(z_1,\dots,z_N)
\end{aligned} 
\hspace{2cm}\text{where}\ z_i=(x_i,y_i),\ i=1,\dots,N,
\right.
\end{equation}
arises naturally as a limit of
weak solutions of the Euler equations for the motion of the whole fluid, see e.g. \cite{flucher-gustafsson}.
The geometry of the domain comes into play through the hydrodynamic Green's function, a
generalisation of the classical Green's function of the first kind for the Laplacian on $\Omega$, which
plays a dominant role in the Hamilton function $H_\Omega$.

Since its derivation by Helmholtz, Kirchhoff, Lord Kelvin and Routh in the second half of the 19\textsuperscript{th} century,
this model has played a role in the research on fluid dynamics, for on the one hand its solutions provide some intuition for the
more general problem of vorticity solutions for the Euler equations, on the other hand its applicability in turbulences of the earth's atmosphere and oceans up to the dynamics of an electron plasma,
see for example the survey article \cite{arefetal} as well as the monographies \cites{marchioro-pulvirenti,newton,saffman}.

There is plenty of literature in the case that $\Omega$ is the whole Euclidean plane (in which case there are only relative equilibria of the system) or all the
vorticities have the same sign. For research about these cases \cites{arefetal,barry-etal,marchioro-pulvirenti,newton,roberts,saffman} 
are a very good starting point. 
In particular, a lot of research has been done concerning existence, stability and geometrical form of stationary or periodic solutions, see for example \cites{arefetal,barry-etal,bartsch-dai,dai-diss,newton2,roberts}.

In these papers it is crucial that if $\Omega$ is the whole Euclidean plane, the Hamilton function is explicitly given by
\[
H(z_1,\dots,z_N)
 =-\frac{1}{2\pi}\sum_{\substack{i,j=1\\ i\ne j}}^N
   \Gamma_i\Gamma_j\ln|z_i-z_j|.
\]
Much less is known about the behaviour of solutions of \eqref{1}, if $\Omega$ is a bounded domain. 
$H_\Omega$ is then defined on the so-called configuration space 
\[
\CF_N\Omega:=\{(z_1,\dots,z_N)\in\Omega^N: z_i\not=z_j \ \text{for}\ i\not=j\},
\]
which is an open subset of $\Omega^N$ and therefore of all of $\IC^{N}$.

In case $\Omega$ is not the whole Euclidean plane, the Hamilton function, which in the literature is commonly called ``Kirch\-hoff--Routh path function'' is given by
\[
H_\Omega(z_1,\dots,z_N)
 = \sum_{j=1}^N \Gamma_j^2h(z_j)
    + \sum_{\substack{i,j=1\\ i\not=j}}^N\Gamma_i\Gamma_j G(z_i,z_j),
\]
where $G:\CF_2\Omega\to\IR$ is the hydrodynamic Green's function with regular part $g(x,y)=G(x,y)+\tfrac{1}{2\pi}\ln|x-y|$, and $h(x)=g(x,x)$
is the so-called Robin's function. 

Basic results concerning $G$, $h$ and the dynamics of the equation \eqref{1} 
in the case $N=1$ were proven in \cites{flucher,flucher-gustafsson,gustafsson}.

Through the Robin's function $h$, interactions of the vortices with the domain's boundary are introduced,
which sets the problem covered here apart from the research on the case that $\Omega=\IC$, where for instance absolute equilibria are inexistent and instead one
is interested in relative equilibria, that is, solutions, where the distance between all vortices remains fixed over time.
If, instead, $\Omega$ is a bounded domain, the presence of $h$ allows for \emph{absolute} equilibria of equation \eqref{1}, which is the topic covered here.

For $N\ge2$ there exist papers of numerical nature by mathematicians, physicists and engineers, especially for special domains, whose Green's function is 
either explicitly known or can be described by methods of complex analysis, we refer to chapter three of \cite{newton} and the references therein.

Contrary to that, there are only a few analytical papers concerned with the case of a general domain $\Omega$, notably \cites{bartolucci-pistoia:2007,ba-pi-we,ba-pi-we-err,delpino-etal:2005,esposito-etal:2005}, where, under some special conditions on $\Omega$ and the coefficients $\Gamma_i$, critical points of $H_\Omega$ and thereby stationary
solutions of \eqref{1} are obtained. In \cite{delpino-etal:2005} it is assumed that $\Gamma_i=1$ for all $i\in\{1,\dots,N\}$ and that $\Omega$ 
is not simply connected. In \cite{esposito-etal:2005} very special simply connected (``dumbbell shaped'') domains are allowed, but again only for the case $\Gamma_i=1$ for all $i\in\{1,\dots,N\}$.

Even less is known if some of the vorticities $\Gamma_i$ are positive, others negative, so that the vortices are rotating in different
directions. This is due to the fact that the Hamiltonian $H_\Omega$ becomes indefinite in this case, in which case the methods of the previously mentioned papers don't apply.

The paper \cite{bartolucci-pistoia:2007} is concerned with the case $N=2$, $\Gamma_1=-1$, $\Gamma_2=1$, and $\Omega$ an arbitrary bounded domain, 
this is the first instance where a stationary configuration of counterrotating vortices in an arbitrary domain is found. 
Lastly, in \cite{ba-pi-we} a stationary solution of $N$ counterrotating vortices lying on the symmetry axis of a axially symmetric domain is found for arbitrary $N$ and $\Gamma_i=(-1)^i$ for $i\in\{1,\dots,N\}$. Additionally, the much more complicated case of a general bounded domain $\Omega$ with $\Gamma_i=(-1)^i$ is settled in \cite{ba-pi-we} successfully for $N=3$ and $N=4$ with corrections and extensions to other values of $\Gamma_i$ in \cite{ba-personal}.

It shall also be mentioned that, somewhat surprisingly, in the papers \cites{bartolucci-pistoia:2007,delpino-etal:2005,esposito-etal:2005}
the Hamiltonian $H_\Omega$ appears as a limit functional for some elliptic boundary value problems in $\Omega$ and the existence of critical
points of certain perturbations of $H_\Omega$ gives rise to solutions of these problems.
Moreover, in \cite{cao-liu-wei-1,cao-liu-wei-2} Cao et al. proved that an isolated stable critical point $(z_1,\dots,z_N)$ of $H_\Omega$ leads to a family $v_\epsilon$ of stationary vector fields in $\Omega$ solving the Euler equations for an ideal fluid such that the vorticities concentrate in blobs near $z_1,\dots,z_N$. As $\epsilon\to0$ the vorticity blobs converge towards the stationary point vorticies $z_1,\dots,z_N$. Although the critical points we obtain here are not necessarily isolated, so that the main theorems from \cite{cao-liu-wei-1,cao-liu-wei-2} do not apply, the methods from \cite{cao-liu-wei-1,cao-liu-wei-2} can be applied as in theorem 2.4 of \cite{ba-personal}.

The goal of this article is to investigate the existence and properties of critical points of $H_\Omega$ (and hence of stationary solutions
to equation \eqref{1}) under various conditions on the vorticities $\Gamma_i$ as well as some geometrical and topological assumptions on $\Omega$, but for general $N\in\IN$. 

In particular, we will prove more general versions of the following two theorems.
\begin{thm}\label{general-pre-result1}
Let $\Gamma\in\IR^N$ satisfy $\sum_{\substack{i,j\in J\\ i\not=j}}\Gamma_i\Gamma_j\not=0$ and $\sum_{j\in J}\Gamma_j^2>\sum_{\substack{i,j\in J\\ i\not=j}}|\Gamma_i\Gamma_j|$ for all $J\subset\{1,\dots,N\}$, $|J|\ge 2$, as well as $\Gamma_j=(-1)^j|\Gamma_j|$, $j\in\{1,\dots,N\}$, where $|\Gamma_j|\le|\Gamma_{j+1}|$ for $j\in\{1,\dots,N-1\}$. Then $H_\Omega$ has a critical point.
\end{thm}
If $\Omega$ is not simply--connected, we may exploit the richer topology of the configuration space $\CF_N\Omega$ to abolish the necessity for alternating vorticities.
\begin{thm}\label{general-pre-result2}
Assume $\Omega$ is not simply--connected and let, as before, $\Gamma\in\IR^N$ satisfy $\sum_{\substack{i,j\in J\\ i\not=j}}\Gamma_i\Gamma_j\not=0$ and $\sum_{j\in J}\Gamma_j^2>\sum_{\substack{i,j\in J\\ i\not=j}}|\Gamma_i\Gamma_j|$ for all $J\subset\{1,\dots,N\}$, $|J|\ge 2$. Then $H_\Omega$ has a critical point.
\end{thm}
More general versions of the above theorems are found in theorems \ref{general} and \ref{notsimplyconnected}, respectively.

The proofs of these theorems seem to suggest that the critical point of theorem \ref{general-pre-result1} has the vortices $z_1,\dots,z_N$ aligned on a line through
domain $\Omega$ such that the signs of the vorticities are alternating along the line, in analogy to the so-called Mallier--Maslowe row of counterrotating vortices (see \cite{mallier-maslowe}),
and that the vortices in the critical point of theorem \ref{general-pre-result2} are placed around a hole in $\Omega$, such that the attraction/repulsion of neighboring
vortices cancel each other out together with the repulsion of the boundary, but unfortunately we are able to prove neither of these results.

However, if the domain has some symmetries, one is able to prove the above conjectures, see for example \cites{ba-pi-we,sym_part}

Another difficult and interesting problem is the question of stability of the derived critical points. As we shall see in the proof, we are embedding $N-2$-- (or $N$ in the case of theorem \ref{general-pre-result2}) dimensional submanifolds into $\CF_N\Omega$, which is a $2N$--dimensional manifold, and deforming them along the gradient flow of $H_\Omega$, therefore the derived critical points should have Morse indices $N+2$ and $N$, respectively and therefore be unstable, provided they are nondegenerate, which in turn is a difficult problem on its own and is not covered here.


Although the problem of finding critical points of $H_\Omega$ is finite dimensional, the problem has proven itself to be considerably refractory.
The most obvious difficulty is that for an arbitrary domain $\Omega$ the Green's function as an essential part of $H_\Omega$ is only implicitly given as
a solution of a partial differential equation, thus all relevant properties of $H_\Omega$ have to be derived through the analysis of the corresponding partial
differential equation.
More importantly, $H_\Omega$ is only defined on the incomplete manifold $\CF_N\Omega$ and is for general vorticities $\Gamma_i$ strongly indefinite. In fact, it
may be the case that $H_\Omega(z)$ remains bounded for
$\dist(z,\partial\CF_N\Omega)\to 0$, which in the model corresponds to collisions of multiple vortices with each other or with $\partial\Omega$. 
This lack of compactness is crucial, since it prevents us from using standard methods of critical point theory, such as the ``mountain--pass''--theorem.
Hence a more detailed study of the behaviour of $H_\Omega$ is necessary. The usual methods of critical point theory, all of which apply some sort of modified
gradient flow of $H_\Omega$ are difficult to apply due to the incompleteness of $\CF_N\Omega$. Success in applying these methods is therefore intimately connected 
to a good analytical understanding of collisions, that is of flow lines $z:(t^-(z_0),t^+(z_0))\to\CF_N\Omega$ to the gradient flow of $H_\Omega$ satisfying
\[\min\left\{|z_i(t)-z_j(t)|,\dist(z_i(t),\partial\Omega):i,j\in\{1,\dots,N\},i\not=j\right\}\to 0\] for $t\to t^+(z_0)$.
This, in turn, depends very sensitively on the constellation of the vorticities $\Gamma_i$.

The space $\CF_N\Omega$ on the other hand exhibits a rich topology even for simply connected $\Omega$, such that, given appropriate compactness properties of $H_\Omega$, finding critical points of $H_\Omega$ is a rather easy task.

The bulk of this article is therefore concerned with deriving conditions on the vorticities $\Gamma_i$ and on the domain $\Omega$ such that the gradient flow of $H_\Omega$
has a compact flow line. The relevant condition on the $\Gamma_i$ has in part already been conjectured in \cite{ba-pi-we} and is a rather strict one for larger $N$. 
In particular, for general $\Omega$, the "model case" $\Gamma_i=(-1)^i$ is not covered by our results, which therefore complement the results given in \cite{ba-pi-we}.

This paper is organized as follows. In section \ref{results_fd} we lay out suitable definitions to concisely state our main results, which are found in subsection \ref{results}. In \ref{preliminaries} we prove some preliminary results concerning the behaviour of the Green's and Robin--functions. We also
give an abstract deformation argument which will be perpetually used throughout the whole article.
Section \ref{closelookchapter} is concerned with the careful analysis of the behaviour of $H_\Omega$ along ``collision'' flowlines. Section \ref{linking_arguments} then provides linking properties for $H_\Omega$, consequently proving the existence of critical points of $H_\Omega$ also in the case of a general domain $\Omega$. 

As far as the general topological techniques of critical point theory such as ``linking phenomena'' are concerned, the book \cites{schechter} provides a good introduction.

The results in this paper were derived in the author's doctoral dissertation \cite{diss}, supervised by T. Bartsch, Giessen, to whom the author wishes to express his thankfulness.

\section{Statement of results}\label{results_fd}
In this section we give an outline of the theorems proven in this paper. 
We start out by fixing some notation and collecting hypotheses on $\Omega$ and $G$. We are then able to concisely state our results in subsection \ref{results}.

\subsection{Hypotheses and basic notation}
\begin{assumption}\label{assumptiononOmega}
Let $\Omega\subset\IC$ be a bounded domain with $C^3$--boundary. A fortiori, $\Omega$ is  finitely connected and satisfies a uniform exterior ball condition, that is there exists a constant $r>0$ such that for any $x\in\partial\Omega$ there is
$x^*\in\IC$ such that $B_r(x^*)\subset\IC\setminus\Omega$ as well as $x\in\partial B_r(x^*)$.
Let $k_0:=\rank\pi_1(\Omega)$, and denote the bounded components (if any) of $\IC\setminus\bar\Omega$ by $\Omega_j$, $j\in\{1,\dots,k_0\}$.
\end{assumption}

For convenience in stating our results and further hypotheses, we start by fixing some useful notation.
\begin{definition}[Basic notation]\label{linesnotation}
The configuration space of $N$ point vortices in $\Omega$ is defined as\[
\CF_N\Omega:=\left\{z\in\Omega^N:z_i=z_j\gdw i=j\right\},
\]
which is an open subset of $\Omega^N$ and therefore of all of $\IC^N$. We denote its boundary in $\IC^N$ by $\partial\CF_N\Omega$.
We set $\Delta_N:=\left\{t\in(0,\infty)^N:t_i<t_{i+1}\ \fa i\in\{1,\dots,N-1\}\right\}$, and for $a\in\IC$ and $v\in S^1$ we define the space of ordered configurations of $N$ vortices along the line
$a+\IR\cdot v$ through $\Omega$ to be the $N$--dimensional submanifold of $\CF_N\Omega$ defined by\[
\CL_N(a,v):=(\tilde a+\Delta_N\cdot v)\cap\Omega^N=\left\{(a+t_1v,\dots,a+t_nv)\in\Omega^N:(t_1,\dots,t_n)\in\Delta_N\right\},
\]
where $\tilde{a}=(a,\dots,a)\in\IC^N$.
The symmetric group $\Sigma_N$ on $N$ symbols acts freely on $\CF_N\Omega$ via\[
\Sigma_N\times\CF_N\Omega\ni(\sigma,z)\mapsto\sigma\ast z:=(z_{\sigma^{-1}(1)},\dots,z_{\sigma^{-1}(N)})\in\CF_N\Omega,
\]
hence it is possible to define\[
\CL_N^\sigma(a,v):=\sigma^{-1}\ast\CL_N(a,v)
\]
as well as
\[
\CL_N^\sigma\Omega:=\bigcup_{(a,v)\in\Omega\times S^1}\CL_N^\sigma(a,v)
\]
for $\sigma\in\Sigma_N$. Lastly, for $\emptyset\not=C\subset\{1,\dots,N\}$ it will be useful to define the
orthogonal projection $\pi_C$ by\[
\pi_C:\IC^N\ni z\mapsto (z_j)_{j\in C}\in\IC^{|C|}.
\]
\end{definition}
\begin{figure}[ht]
\begin{center}
\includegraphics[scale=1]{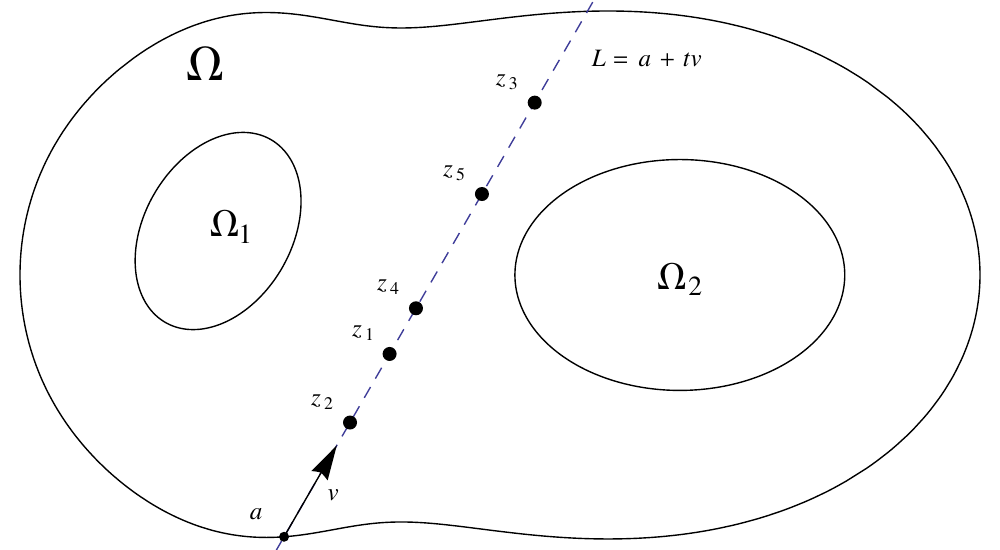}
\end{center}
\caption{A configuration $z\in\CL_5^{\sigma}(a,v)\subset\CF_5\Omega$, where $\sigma=(12)(345)$, compare def. \ref{linesnotation}.}
\end{figure}

\begin{definition}[Reflection at $\partial\Omega$]
Since $\Omega$ is $C^3$ there is $\eps>0$ such that the orthogonal projection\[
p:\Omega_\eps:=\{z\in\Omega:\dist(z,\partial\Omega)<\eps\}\to\partial\Omega
\]
is a well--defined $C^2$--map satisfying $|p(z)-z|=\dist(z,\partial\Omega)$. The reflection at $\partial\Omega$ is then defined
as the $C^2$--map
\[
\overline{\Omega_\eps}\ni z\mapsto\overline{z}:=2p(z)-z\in\IC.
\]
Here and in all what follows, when talking about differentiability, we regard $\IC=\IR^2$, and thus mean differentiability in the real--valued sense. We regard $\Omega\subset\IC$ simply because the elegant geometrical properties of complex multiplication allow us to state some things more concisely.
Additionally, in what comes we will always abbreviate $d(z):=\dist(z,\partial\Omega)$.
\end{definition}

We are now ready to state the general sufficient assumptions on the function $G$ for carrying out our arguments.

\begin{assumption}\label{assumptiononG}
Let $G:\CF_2\Omega\to\IR$
 satisfy the following hypotheses: $G$ is bounded below by some constant $C_0$ and has logarithmic singularities
on the diagonal in $\Omega\times\Omega$, more precisely,
the map $\CF_2\Omega\ni(x,y)\mapsto G(x,y)+\frac{1}{2\pi}\ln|x-y|\in\IR$ has a continuation $g\in C^1(\Omega^2)$, which is bounded from above by some constant $C_1>0$.
Thus, we may write
\begin{equation}\label{log-singularity}
G(x,y)=g(x,y)-\frac{1}{2\pi}\ln|x-y|.
\end{equation}
Further, for every $\eps>0$ there is $C_2>0$ depending only on $\Omega$ and $\eps$ such that
\begin{equation}\label{C1bound}
|G(x,y)|+|\nabla_x G(x,y)|+|\nabla_y G(x,y)|\le C_2
\end{equation}
for every $x,y\in\Omega$ with $|x-y|\ge\eps$. Similarly, there is a constant $C_3>0$, also depending only on $\eps$ and $\Omega$, such that
\begin{equation}\label{reflection}
|\psi(x,y)|+|\nabla_x \psi(x,y)|+|\nabla_y \psi(x,y)|\le C_3, 
\end{equation}
for every $x,y\in\Omega_\eps$,
where $\psi(x,y)=g(x,y)-\tfrac{1}{2\pi}\ln|\bar x-y|$ and $x\mapsto\bar x$ is reflection at $\partial\Omega$.
Further there exists a constant $C_4>0$ such that for any line $L=\IR v+w\subset\IC$ with $L\cap\Omega\not=\emptyset$ 
\begin{equation}\label{bddonlines_hyp}
G(w+rv,w+sv)-G(w+rv,w+tv)\ge -C_4.
\end{equation}
for all $r<s<t$, for which the left hand side is defined. 
\end{assumption}

Concerning interesting candidates for a function $G$ satisfying the above hypotheses, we have the following result.
\begin{thm}\label{dirichlet}
Green's function of the first kind for the Dirichlet Laplacian in $\Omega$ satisfies hypothesis \ref{assumptiononG}.
\end{thm}

Rather than the regular Green's function for the Dirichlet Laplacian, the single most important class of Green's functions $G$ for fluid dynamics is the class of so-called \emph{hydrodynamic Green's functions}, which we will now introduce. An excellent motivation and introduction to the topic of hydrodynamic Green's functions is provided by \cite{flucher-gustafsson}.

\begin{definition}[Hydrodynamic Green's function]\label{hgf}
The \emph{hydrodynamic Green's function with periods} $\gamma_0,\dots,\gamma_{k_0}\in\IR$, subjected to the condition $\sum_{j=0}^{k_0}\gamma_j=-1$ is the unique solution $G\in C^2(\CF_2\overline{\Omega})$ of the problem 
\[
\begin{cases}
-\Delta G(\cdot,y)=\delta_y & \text{for every}\ y\in\Omega \\
\spr{\nabla_xG(x,y),\tau_x}=0 & \text{for every}\ y\in\Omega, x\in\partial\Omega \\
\int_{\partial\Omega_j}\spr{\nabla_xG(x,y),\nu_x}\dif s(x)=\gamma_j & \text{for every}\ j\in\{0,\dots,k_0\}\\
\int_{\partial\Omega}G(x,y)\spr{\nabla_xG(x,z),\nu_x}\dif s(x)=0 & \text{for every}\ y,z\in\Omega, 
\end{cases}
\]
where $\partial\Omega_0=\partial\Omega\setminus\bigcup_{j=1}^{k_0}\partial\Omega_j$.
\end{definition}
Using this definition we have the following 
\begin{thm}\label{hydrodynamic_gf}
Any hydrodynamic Green's function satisfies hypothesis \ref{assumptiononG}.
\end{thm}

We postpone the proof of the last two theorems until the next section. Having these results in mind, we will in the following by slight abuse of language refer to any function $G$ satisfying hypothesis \ref{assumptiononG} as a (hydrodynamic) Green's function on $\Omega$.

\begin{definition}
For $\Gamma\in\IR^N$ we define the \emph{Kirchhoff--Routh path function} for vortices with vorticities $\Gamma_i$, $i\in\{1,\dots,N\}$ to be the function\[
H_\Omega^\Gamma:\CF_N\Omega\ni(z_1,\dots,z_N)\mapsto\sum_{j=1}^N\Gamma_j^ 2h(z_j)+\sum_{\substack{i,j=1\\ i\not=j}}^N\Gamma_i\Gamma_jG(z_i,z_j)\in\IR,
\]
where the function $G$ satisfies hypothesis \ref{assumptiononG} and $h(x)=g(x,x)$ for all $x\in\Omega$. If the parameter $\Gamma$ is understood, we will drop it from the notation, writing $H_\Omega$ instead of 
$H_\Omega^\Gamma$.
\end{definition}

Let us conclude this subsection by specifying the necessary technical conditions on the parameter $\Gamma$.

\begin{definition}[$\Delta$--admissibility]\label{delta}
We call a parameter $\Gamma\in\IR^N$ $\Delta$--admissible, if for every $C\subset\{1,\dots,N\}$, $|C|\ge 2$: \[\sum\limits_{\substack{i,j\in C\\i\not=j}}\Gamma_i\Gamma_j\not=0.\]
\end{definition}
\begin{definition}[$\partial$--admissibility]\label{partial}
We call a parameter $\Gamma\in\IR^N$ $\partial$--admissible, if for every $C\subset\{1,\dots,N\}$, $|C|\ge 2$: \[
\sum_{i\in C}\Gamma_i^2>\sum_{\substack{i,j\in C\\ i\not=j}}|\Gamma_i\Gamma_j|.\]
If $\Omega$ is strictly convex, this condition may be replaced by\[
\sum_{i\in C}\Gamma_i^2>\sum_{\substack{i,j\in C\\ \Gamma_i\Gamma_j<0}}|\Gamma_i\Gamma_j|.
\]
\end{definition}

Without any geometric or topological assumptions on $\Omega$ we need more specialised parameters $\Gamma$ in order to derive positive results. These conditions are stated in the following definition, which shall conclude this subsection.
\begin{definition}[$\CL$--admissible parameters]\label{L-admissible}
A parameter $\Gamma\in\IR^N$ is called $\CL$-admissible if there is $\sigma\in\Sigma_N$ such that $\iota(\sigma\ast\Gamma)\in\overline{\Delta_N}$ or $-\iota(\sigma\ast\Gamma)\in\overline{\Delta_N}$, where 
$\iota:\IR^N\to\IR^N$ is the involution $(x_j)_{j\in\{1,\dots,N\}}\mapsto\left((-1)^{j}x_j\right)_{j\in\{1,\dots,N\}}$ and the closure is to be taken in $(0,\infty)^N$. Similarly, we call $\Gamma$ strictly
$\CL$--admissible, if $\iota(\sigma\ast\Gamma)\in\Delta_N$ or $-\iota(\sigma\ast\Gamma)\in\Delta_N$.
\end{definition}

The intuition behind the definitions \ref{delta} and \ref{partial} is simple: the idea is that $\Delta$--admissibility prevents collisions of vortices inside the ``diagonal'' $\Delta=\{z\in\Omega^N:z_i=z_j\ \text{for some}\ i\not=j\}$ from happening while the energy $H_\Omega$ of the system remains finite. Since $\Gamma_i\Gamma_jG(z_i,z_j)$
becomes large if $z_i$ and $z_j$ collide inside $\Omega$, we may regard the quantity $\sum_{\substack{i,j\in C\\i\not=j}}\Gamma_i\Gamma_j$ as a kind of ``collision weight''
associated to the vortices with indices in $C$. Since $\Gamma_j^2 h(z_j)\to -\infty$ if $z_j\to\partial\Omega$, we may, by the same intuitive reasoning, regard the quantity
$\sum_{i\in C}\Gamma_i^2$ as a kind of weight for the boundary interaction of the vortices $z_i$, $i\in C$, and the condition of $\partial$--admissibility then simply states that the boundary interaction outweighs the collision weight, a condition which presents itself naturally in the proof of lemma \ref{nocollisionwithbdry}, as we may see later on.

The condition of $\CL$--admissibility means, that one is able to align the vortices $z_1,\dots,z_N$ along a line in such a way that the signs of the vorticities are alternating and their moduli nonincreasing. The motivation behind this is the physical intuition that it should be possible for vortices aligned on an axis, that the repulsion between two neighbouring vortices of opposite sign is cancelled out by the attraction to the next but one vortices, whose vorticities have the same sign again, in analogy to the well-known Mallier--Maslowe row of counterrotating vortices.

\subsection{Main results}\label{results}
Equipped with these definitions, we are now able to concisely state the theorems proven in this paper. 

\begin{thm}\label{general}
For any $N\in\IN$ and any $\CL$--admissible, $\partial$--admissible and $\Delta$--admissible parameter $\Gamma\in\IR^N$, the Kirchhoff--Routh path function $H_\Omega$ has a critical point in $\CF_N\Omega$.
\end{thm}

If $\Omega$ is not simply--connected, we may exploit the richer topology of the configuration space $\CF_N\Omega$ to abolish the necessity for alternating vorticities.

\begin{thm}\label{notsimplyconnected}
Suppose that $\pi_1(\Omega)\not=0$. Then for any $N\in\IN$ and for any $\partial$--admissible and $\Delta$--admissible parameter $\Gamma\in\IR^N$, the Kirchhoff--Routh path function $H_\Omega$ has a critical point in $\CF_N\Omega$.
\end{thm}

As work on these results started out, the paper \cite{ba-pi-we} was the main starting point for the research conducted here. As of this writing, \cite{ba-pi-we} and the recent preprint \cite{ba-personal} are the only references known to the author dealing with $N\ge3$ vortices having both general and alternating signs of vorticities. Without symmetries of the domain only the cases $N=3$ and $N=4$ have been treated in \cite{ba-pi-we,ba-personal}. For larger values of $N$ the problem turns out to be a much more difficult one. The paper at hand  is a first step in this direction providing criteria on the vorticity vector $\Gamma$ which guarantee the existence of critical points of $H_\Omega$ for arbitrary $N$. The author believes, partly based on numerical simulations of the problem, that the severe restriction of $\partial$-admissibility is in fact unnecessary and may be completely abolished or at least be weakened, as it is done in \cite{ba-pi-we,ba-personal} for $N\le4$. A detailed investigation of the influence of symmetries going far beyond the results in \cite{ba-pi-we} can be found in \cite{diss,sym_part}.

\section{Preliminaries}\label{preliminaries}

\subsection{Preliminary results}
This section is concerned with the proof of some preliminary results. We start out by proving the preliminary theorems \ref{dirichlet} and \ref{hydrodynamic_gf}.

\begin{proof}[\textsc{Proof} (of theorem \ref{dirichlet})]
All of the conditions in \ref{assumptiononG} are either well known properties of the Dirichlet Laplacian, see for example \cite{gilbarg-trudinger} or verified in \cite{ba-pi-we} except for property \eqref{bddonlines_hyp}, which is a slight sharpening of the result given there. 

To see that \eqref{bddonlines_hyp} holds assume on the contrary that there is a sequence $L_n=a_n+\IR v_n$ of lines with $\Omega\cap L_n\not=\emptyset$ as well as $r_n<s_n<t_n$ such that\begin{equation}\label{contr}
G(a_n+r_nv_n,a_n+s_nv_n)-G(a_n+r_nv_n,a_n+t_nv_n)\to -\infty
\end{equation}
as $n\to\infty$, where by selecting appropriate subsequences we may take the sequences $(a_n)\subset\overline{\Omega}$ and $(v_n)\subset S^1$ to be convergent to some $a\in\overline{\Omega}$ and $v\in S^1$, respectively.

Now since $G$ is bounded below \eqref{contr} implies
\[
G(x_n,z_n)=g(x_n,z_n)-\frac{1}{2\pi}\ln|t_n-r_n|\to\infty,
\]
hence $|t_n-r_n|\to 0$, such that $r_n,s_n,t_n\to t$ as $n\to\infty$, since $g$ is bounded from above, and where we abbreviated $x_n:=a_n+r_nv_n$, $y_n:=a_n+s_nv_n$, $z_n:=a_n+t_nv_n$.

If $a+tv\in\Omega$ this leads to a contradiction via\[
G(x_n,y_n)-G(x_n,z_n)\ge C-\frac{1}{2\pi}\ln\frac{s_n-r_n}{t_n-r_n}\ge C-\ln 1=C,
\]
for some constant $C$, since $g$ is bounded on compact subsets of $\Omega\times\Omega$. 

Thus $a+tv\in\partial\Omega$, so if $n$ is large enough, $x_n,y_n,z_n\in\Omega_\eps$, where we may use the approximation\[
G(x_n,y_n)-G(x_n,z_n)=\frac{1}{2\pi}\ln\frac{|\overline{x_n}-y_n|}{|\overline{x_n}-z_n|}\cdot\frac{|x_n-z_n|}{|x_n-y_n|}+O(1)
\]
as $n\to\infty$. Considering the differentiable function\[
f:(r_n,\infty)\ni\alpha\mapsto\frac{|\overline{x_n}-a_n-\alpha v_n|^2}{|x_n-a_n-\alpha v_n|^2}\in\IR
\]
we easily compute\[
f'(\alpha)=2\frac{\spr{-v_n,\overline{x_n}-a_n-\alpha v_n}|x_n-a_n-\alpha v_n|^2+\spr{v_n,x_n-a_n-\alpha v_n}|\overline{x_n}-a_n-\alpha v_n|^2}{|x_n-a_n-\alpha v_n|^4}
\]\[
=\frac{(4d_{x_n}\spr{v_n,\nu_{x_n}}+\alpha-r_n)(\alpha-r_n)^2+2\spr{v_n,r_nv_n-\alpha v_n}|\overline{x_n}-a_n-\alpha v_n|^2}{|x_n-a_n-\alpha v_n|^4}\]\[\begin{split}
=\frac{2}{|x_n-a_n-\alpha v_n|^4}\bigg(
(2d_{x_n}\spr{v_n,\nu_{x_n}}(\alpha-r_n)^2+(\alpha-r_n)^3\\
+(r_n-\alpha)\left(4d_{x_n}^2+(r_n-\alpha)^2-2d_{x_n}(r_n-\alpha)\spr{v_n,\nu_{x_n}}\right)
\bigg)\end{split}
\]\[
=\frac{8d_{x_n}^2(r_n-\alpha)}{|x_n-a_n-\alpha v_n|^4}\le 0,
\]
thus $f$ is decreasing, in other words\[
\frac{|\overline{x_n}-y_n|^2}{|x_n-y_n|^2}=f(s_n)\ge f(t_n)=\frac{|\overline{x_n}-z_n|^2}{|x_n-z_n|^2},
\]
hence\[
\frac{|\overline{x_n}-y_n|}{|\overline{x_n}-z_n|}\cdot\frac{|x_n-z_n|}{|x_n-y_n|}\ge 1,
\]
from which we deduce\[
G(x_n,y_n)-G(x_n,z_n)=\frac{1}{2\pi}\ln\frac{|\overline{x_n}-y_n|}{|\overline{x_n}-z_n|}\cdot\frac{|x_n-z_n|}{|x_n-y_n|}+O(1)
\]
\[
\ge \ln 1+O(1)=O(1),
\]
as $n\to\infty$, which is the desired contradiction.
\end{proof}

We continue by proving theorem \ref{hydrodynamic_gf}.

\begin{proof}[\textsc{Proof} (of theorem \ref{hydrodynamic_gf})]
Nearly all of this follows from the fact that there is a symmetric positive semidefinite matrix $(g^{kl})\in\IR^{k_0\times k_0}$, such that\[
G(x,y)=G^0(x,y)+\sum_{k,l=1}^{k_0+1}g^{kl}u_{k-1}(x)u_{l-1}(y),
\]
where $G^0$ is the Green's function of the Dirichlet Laplacian in $\Omega$ and the $u_k$ are the unique solutions of\[
\begin{cases}
\Delta u_k=0 & \text{in}\ \Omega\\
u_k=\delta_{kl} & \text{on}\ \partial\Omega_l,
\end{cases}
\]
see \cite{flucher-gustafsson}, proposition 7.
By assumption \ref{assumptiononOmega} on $\partial\Omega$ each of the $u_k$ has bounded gradient and is bounded by the maximum principle.
Therefore \eqref{log-singularity}, \eqref{C1bound} and \eqref{bddonlines_hyp} are immediate and so is \eqref{reflection}, since\[
\psi(x,y)=g(x,y)-\tfrac{1}{2\pi}\ln|\bar x-y|=g_0(x,y)-\sum_{k,l}g^{kl}u_k(x)u_l(y)-\tfrac{1}{2\pi}\ln|\bar x-y|,
\]
in other words $\psi(x,y)=\psi_0(x,y)-\sum_{k,l}g^{kl}u_k(x)u_l(y)$ and we are done.
\end{proof}

Concerning the analysis of the boundary behaviour of $H_\Omega$, the condition \eqref{reflection} is of course crucial. The detailed study of boundary collisions will be postponed until chapter \ref{closelookchapter}, but by then we will need a technical lemma, which may be a simple case of some general theorem known to differential geometers.
\begin{lemma}\label{Dp} There is $\eps>0$ such that\[
\Dif p(z)v=\frac{1}{1-\kappa_zd_z}\spr{v,\tau_z}\tau_z
\]\[
\Dif\nu_zv=-\frac{\kappa_z}{1-\kappa_zd_z}\spr{\tau_z,v}\tau_z
\]
holds for any $z\in\Omega_\eps$,
where $\tau_z$ is the unit tangent vector to $\partial\Omega$ at $p(z)$, such that the basis $(\tau_z,\nu_z)$ is positively oriented and $\kappa_z$ is the curvature of $\partial\Omega$ at $p(z)$ with respect to the induced orientation of $\partial\Omega$.
\end{lemma}

\begin{proof2}
This is an easy application of the implicit function theorem.
\end{proof2}

\begin{lemma}\label{distinequalities}
There are $\eps>0$ and constants $C_6,C_7,C_8>0$ depending only on $\Omega$ such that the inequalities\begin{equation}
\max\{d(x)+d(y),C_6|x-y|\}\le|x-\overline{y}|\le |x-y|+2d(y)
\end{equation}
\begin{equation}
|x-\overline{y}|^2\ge C_7|p(x)-p(y)|^2
\end{equation}
\begin{equation}
\left||x-\overline{y}|-|\overline{x}-y|\right|^2\le C_8(d(x)+d(y))|p(x)-p(y)|^2
\end{equation}
hold for any $x,y\in\Omega_\eps$.
\end{lemma}

\begin{proof2}
Concerning the first inequality consider the straight line joining $x$ and $\overline{y}$. This line intersects $\partial\Omega$ at some point $z\in\partial\Omega$, which implies\[
|x-\overline{y}|=|x-z|+|z-\overline{y}|\ge d(x)+d(\overline{y})=d(x)+d(y).
\]
The other direction is immediate from the triangle inequality, since\[
|x-\overline{y}|=|x-y+y-\overline{y}|=|x-y+2d(y)\nu_y|\le |x-y|+2d(y),
\]
and the other inequalities are verified as (2.1), (2.4) and (2.5) in \cite{ba-pi-we}. 
\end{proof2}

With this notation, we have the following theorem, which lies on the very foundation of this thesis. Its proof is similar to the one given for the case of axially symmetric
$\Omega$ in \cite{ba-pi-we} but works out just as well for general $\CL$--admissible parameters $\Gamma$ without any assumptions on symmetry.
\begin{thm}\label{bddonlinesandsoon}
Let $\Gamma$ be $\CL$--admissible with corresponding permutation $\tilde\sigma\in\Sigma_N$ and let
$\sigma\in\{\hat\sigma\tilde{\sigma},\tilde{\sigma}\}$, where $\hat\sigma\in\Sigma_N$ is the order--reversing permutation. Then $H_\Omega\big|_{\CL_N^\sigma\Omega}$ is bounded above, and fixing a line 
$L=a+\IR v\subset\IC$ with $a\in\IC\setminus\Omega$, $v\in S^1$ and $\Omega\cap L\not=\emptyset$, we have that
\[H_\Omega\big|_{\CL_N^\sigma(a,v)}(z)\to-\infty \quad\text{as}\quad z\to\partial\CL_N^\sigma(a,v),\]
where the boundary of the $N$--dimensional submanifold $\CL_N^\sigma(a,v)$ of $\CF_N\Omega$ is to be taken in $L^N$. 
\end{thm}

\begin{proof2}
Let $\Gamma$ be $\CL$--admissible. Since the change $\Gamma\mapsto-\Gamma$ leaves 
$H_\Omega$ unaffected we may assume without loss of generality that $\Gamma_j=(-1)^j|\Gamma_j|$ and $|\Gamma_{j+1}|\le|\Gamma_j|$. 
Thus $H_\Omega$ takes the form\[
H_\Omega(x)=\sum\limits_{j=1}^N \Gamma_j^2h(x_j)+\sum\limits_{\substack{i,j=1\\ i\not=j}}^N(-1)^{i+j}|\Gamma_i\Gamma_j|G(x_i,x_j)
\]\[
=\sum\limits_{j=1}^N \Gamma_j^2h(x_j)+2\sum\limits_{i=1}^{N-1}G_i(x),
\]
where for $N-i$ even we have
\[
G_i(x)=\sum\limits_{k=1}^{\frac{N-i}{2}}|\Gamma_i|\bigg(|\Gamma_{i+2k}|G(x_i,x_{i+2k})-|\Gamma_{i+2k-1}|G(x_i,x_{i+2k-1})\bigg),
\]
whereas for $N-i$ odd\[
G_i(x)=\sum\limits_{k=1}^{\frac{N-i-1}{2}}|\Gamma_i|\bigg(|\Gamma_{i+2k}|G(x_i,x_{i+2k})-|\Gamma_{i+2k-1}|G(x_i,x_{i+2k-1})\bigg)-|\Gamma_i\Gamma_N|G(x_i,x_N)
\]
We now are able to infer from hypothesis \eqref{bddonlines_hyp} that for any line $L = \{a + tv : t \in\IR\}$ with $a\in\partial\Omega$, $v\in S^ 1$, and such that $\Omega\cap L\not=\emptyset$\[
G(a + rv, a + sv) - G(a + rv, a + tv) \ge -C_4\]
for all $r < s < t$ for which the left hand side is defined, so combining this result with the condition $|\Gamma_{i-1}|\ge|\Gamma_i|$ and the fact that $G\ge C_0$ we get for a $x\in\CL_N(a,v)$
and $N-i$ even, that $G_i(x)$ is equal to
\[
\sum\limits_{k=1}^{\frac{N-i}{2}}\!|\Gamma_i|\Bigg[|\Gamma_{i+2k}|\bigg(\!\underbrace{G(x_i,x_{i+2k})-G(x_i,x_{i+2k-1})}_{\le C_4}\bigg)
+(\underbrace{|\Gamma_{i+2k}|-|\Gamma_{i+2k-1}|}_{\le 0})G(x_i,x_{i+2k-1})\Bigg]
\]\[
\le \sum\limits_{k=1}^{\frac{N-i}{2}}\bigg(|\Gamma_i\Gamma_{i+2k}|C_4+|\Gamma_i|\left(|\Gamma_{i+2k}|-|\Gamma_{i+2k-1}|\right)C_0\bigg),
\]
whereas analogously for $N-i$ odd\[
G_i(x)\le 
\sum\limits_{k=1}^{\frac{N-i-1}{2}}\bigg(|\Gamma_i\Gamma_{i+2k}|C_4+|\Gamma_i|\left(|\Gamma_{i+2k}|-|\Gamma_{i+2k-1}|\right)C_0\bigg)-\underbrace{|\Gamma_i\Gamma_N|G(x_i,x_N)}_{\ge |\Gamma_i\Gamma_N|C_0}.\]
Since by hypothesis \ref{assumptiononG} $h$ is bounded from above by $C_1$, this gives the required upper bound.

Now fixing $a$ and $v$, every $z\in\CL_N(a,v)$ has a unique representation $z=\tilde a+tv$ with $t\in\Delta_N$, where $\tilde{a}:=(a,\dots,a)\in\IC^N$. Setting\[
\CR:=\left\{t\in\Delta^N:\tilde a+tv\in\CF_N\Omega\right\}
\]
as well as\[
E:\CR\ni t\mapsto H_\Omega(\tilde a+tv)\in\IR,
\]
we have to show that $E(t)\to-\infty$ as $\dist(t,\partial\CR)\to 0$. Therefore consider a sequence $(t^n)\subset\CR$ with the property that
$t^n\to\partial\CR$ as $n\to\infty$. Let us first consider the case that $d(a+t_k^nv)\to 0$ as $n\to\infty$ for some $k\in\{1,\dots,N\}$.
Since $\sum_{i,j=1, i\not=j}^N(-1)^{i+j}|\Gamma_i\Gamma_j|G(a+t_i^nv,a+t_j^nv)$ is bounded from above as $n\to\infty$ and $h(a+t_k^nv)\to -\infty$ for
$n\to\infty$ we infer that indeed $E(t^n)\to-\infty$ as claimed.

Hence we may assume that \begin{equation}\label{bdry}
\liminf_{n\to\infty}d(a+t_j^nv)>0
\end{equation} for all $j\in\{1,\dots,N\}$ and that
\begin{equation}\label{coll}
t_{k+1}^n-t_k^n\to 0 \quad\text{as}\ n\to\infty
\end{equation}
for some $k\in\{1,\dots,N-1\}$. By assumption \eqref{bdry} the first two sums in\[
E(t^n)=\sum_{j=1}^N\Gamma_j^2h(a+t_j^nv)+\sum\limits_{\substack{i,j=1\\ i\not=j}}^N\Gamma_i\Gamma_jg(a+t_i^nv,a+t_j^nv)
-\sum\limits_{\substack{i,j=1\\ i\not=j}}^N(-1)^{i+j}\frac{|\Gamma_i\Gamma_j|}{2\pi}\ln|t_i^n-t_j^n|
\]
remain bounded as $n\to\infty$. We then expand\[
\sum\limits_{\substack{i,j=1\\ i\not=j}}^N(-1)^{i+j}|\Gamma_i\Gamma_j|\ln|t_i^n-t_j^n|
=2\sum_{i=1}^{N-1}\sum_{j=i+1}^N(-1)^{i+j}|\Gamma_i\Gamma_j|\ln|t_i^n-t_j^n|
=2\sum_{i=1}^{N-1}|\Gamma_i|\ln\psi_i(t),
\]
where for $N-i$ even\[
\psi_i(t)=\prod_{j=1}^{\frac{N-i}{2}}\frac{|t_{i+2j}^n-t_i^n|^{|\Gamma_{i+2j}|}}{|t_{i+2j-1}^n-t_i^n|^{|\Gamma_{i+2j-1}|}}
\ge \prod_{j=1}^{\frac{N-i}{2}}|t_{i+2j-1}^n-t_i^n|^{|\Gamma_{i+2j}|-|\Gamma_{i+2j-1}|}\ge C
\]
and for $N-i$ odd\[
\psi_i(t)=\frac{1}{|t_N^n-t_i^n|^{|\Gamma_N|}}\prod_{j=1}^{\frac{N-i-1}{2}}\frac{|t_{i+2j}^n-t_i^n|^{|\Gamma_{i+2j}|}}{|t_{i+2j-1}^n-t_i^n|^{|\Gamma_{i+2j-1}|}}
\ge \frac{C}{|t_N^n-t_i^n|^{|\Gamma_N|}}
\]
for some constant $C>0$, since $|\Gamma_{i+2k}|-|\Gamma_{i+2k-1}|\le 0$.
It thus remains to show that $\psi_k(t^n)\to\infty$ for $n\to\infty$ for some $k\in\{1,\dots,N-1\}$. 
Let $k$ be maximal satisfying \eqref{coll}. If $k=N-1$ we infer $\psi_{N-1}(t^n)\to\infty$ for $n\to\infty$ and the proof is done. Otherwise there is $\delta>0$ such 
that \begin{equation}\label{nocoll}
\delta\le|t_{j+1}^n-t_j^n|\le\frac{1}{\delta} \quad\text{for all}\ j>k
\end{equation}
and $n$ sufficiently large. Therefore, if $N-k$ is even,\[
\psi_k(t)=\frac{|t_{k+2}^n-t_k^n|^{|\Gamma_{k+2}|}}{|t_{k+1}^n-t_k^n|^{|\Gamma_{k+1}|}}\prod_{j=2}^{\frac{N-k}{2}}\frac{|t_{k+2j}^n-t_k^n|^{|\Gamma_{k+2j}|}}{|t_{k+2j-1}^n-t_k^n|^{|\Gamma_{k+2j-1}|}}
\ge \frac{\tilde C\cdot\delta^{|\Gamma_{k+2}|}}{|t_{k+1}^n-t_k^n|^{|\Gamma_{k+1}|}}\to\infty
\]
as $n\to\infty$ by \eqref{coll} and \eqref{nocoll}, whereas for $N-k$ odd\[
\psi_k(t)=\frac{|t_{k+2}^n-t_k^n|^{|\Gamma_{k+2}|}}{|t_{k+1}^n-t_k^n|^{|\Gamma_{k+1}|}}\cdot\frac{1}{|t_N^n-t_k^n|^{|\Gamma_N|}}\prod_{j=2}^{\frac{N-k}{2}}\frac{|t_{k+2j}^n-t_k^n|^{|\Gamma_{k+2j}|}}{|t_{k+2j-1}^n-t_k^n|^{|\Gamma_{k+2j-1}|}}\]\[
\ge \frac{\tilde C\cdot\delta^{|\Gamma_{k+2}|+1}}{|t_{k+1}^n-t_k^n|^{|\Gamma_{k+1}|}}\to\infty\quad \text{as} \ n\to\infty,
\]
again by \eqref{coll} and \eqref{nocoll} and the proof is done.
\end{proof2}

\subsection{A general deformation argument}
For simplicity in stating our results, we find the following definition useful:

\begin{definition}[$\phi$-complete deformations]
Let $X$ be a topological space and let $\phi$ be a flow on $\CF_N\Omega$. 
We call a family $\beta\subset[\alpha]\in[X,\CF_N\Omega]$ of homotopic maps $\phi$-complete, if for any
$\alpha\in\beta$ and any continuous map $T:\alpha(X)\to[0,\infty)$ such that $T(x)\in[0,t^+(x))$ for all $x\in\alpha(X)$ the map\[
X\ni x\mapsto \phi(T(\alpha(x)),\alpha(x))\in\CF_N\Omega
\]
is in $\beta$.
\end{definition}

Denoting the gradient flow of $H_\Omega$ by\[\phi:\bigcup_{z\in\CF_N\Omega}(t^-(z),t^+(z))\times\{z\}\to\CF_N\Omega,\]
in the sequel, we will frequently use the following
\begin{lemma}[general deformation argument]\label{deformation}
Suppose there is a subset $\CL\subset\CF_N\Omega$, such that $H_\Omega$ is bounded above
on $\CL$ by $\sigma$, that is
\begin{equation}\label{bdd}
\sup H_\Omega\big|_{\CL}=\sigma<\infty. 
\end{equation}
Further let $X$ be a topological space, $\beta\subset [\alpha]\in [X,\CF_N\Omega]$ be $\phi$-complete, and such that for any representative 
$\alpha\in \beta$ the intersection 
\[\alpha(X)\cap \CL \not=\emptyset
\] 
is nonempty. Then, fixing some representative $\alpha_0$, there is $x\in\alpha_0(X)$, 
such that 
\[\lim\limits_{t\to t^+(x)}H_\Omega(\phi(t,x))<\infty.\]
\end{lemma}

\begin{proof2}
Assume not. Then for every $x\in\alpha_0(X)$ there is a minimal $T(x)\in[0,t^+(x))$ such that
\[H_\Omega(\phi(T(x),x))=\sigma+1.\]
Since for each $x\in\alpha_0(X)$ $H_\Omega\circ\psi(\cdot,x)$ is strictly increasing (otherwise we are done), the map\[
T:\alpha_0(X)\ni x\mapsto T(x)\in\IR
\]
is continuous. It follows, that unambiguously defining\[
\alpha_t:X\ni\xi\mapsto\phi\left(tT(\alpha_0(\xi)),\alpha_0(\xi)\right)\in\CF_N\Omega,
\]
where $t\in[0,1]$: $\alpha_0\simeq\alpha_1$ and $\alpha_1\in\beta$, since $\beta$ is $\phi$-complete, hence there is 
$\xi\in X$ with $\alpha_1(\xi)\in \CL$, but 
$H_\Omega(\alpha_1(\xi))=\sigma+1$, in contradiction with \eqref{bdd}, and we are done.
\end{proof2}

\section{Singularities of \texorpdfstring{$H_\Omega$}{the functional}}\label{closelookchapter}
This section is devoted to the study of $H_\Omega$ near collisions with the boundary $\partial\Omega$ or with each other away from the boundary and to give conditions on $\Gamma$ and $\Omega$ which prevent these.

The goal of this section is to prove the following
\begin{prop}\label{nocollisions}
Let $\Gamma\in\IR^N$ be $\partial$--admissible and $\Delta$--admissible. Then there is $\delta>0$ such that
$|\nabla H_\Omega(z)|>1$ for every $z$ in $\CM_\delta$,
where\[
\CM_\delta:=\left\{z\in\CF_N\Omega:|z_i-z_j|\le\delta\ \text{or}\ d(z_j)\le\delta\ \text{for some}\ i,j\in\{1,\dots,N\},i\not=j\right\}.
\]
In particular, $H_\Omega$ satisfies the Palais--Smale--condition. Further, if there is $z\in\CF_N\Omega$ such that\[
\lim\limits_{t\to t^+(z)}H_\Omega(\phi(z,t))<\infty,
\]
we have $t^+(z)=\infty$ and there is a sequence $s_n\to\infty$ such that defining $z^{s_n}:=\phi(z,s_n)$, we have $z^{s_n}\to z^*\in\CF_N\Omega$ where $\nabla H_\Omega(z^*)=0$, hence $H_\Omega$ has a critical point.
\end{prop}
\enlargethispage{1ex}

The proof of proposition \ref{nocollisions} of course involves a detailed study of the behaviour of $H_\Omega$ near its singularities.
The functional $H_\Omega$ has singularities at the boundary $\partial\CF_N\Omega$ of $\CF_N\Omega$ in $\IC^N$. This boundary consists of points $z\in\overline\Omega^N$ with $z_j\in\partial\Omega$ or $z_i=z_j$ for some indices $i,j\in\{1,\dots,N\}$, $i\not=j$, corresponding to collisions of vortices with the boundary or with each other in $\overline{\Omega}$, respectively.

In order to deal with the problem of collisions effectively, we first introduce some convenient notation for dealing with different types of collisions of vortices
within $\overline\Omega$, corresponding to the respective parts of the boundary $\partial\CF_N\Omega$ of $\CF_N\Omega$ in $\IC^N$.In order not to get too deep into technicalities already in the introduction, we state the results in a simplified rather than their fully general version in order to give an overview of the topics covered.

First note that collisions of vortices correspond to partitions of the set $\{1,\dots,N\}$ as follows: 
Given a point $z\in\overline{\Omega}^N$, we define \[\CP_z:=\left\{
C\subset\{1,\dots,N\}: z_i=z_j\gdw i,j\in C
\right\},\]
which is clearly a partition of $\{1,\dots,N\}$.
We call an element $C\in\CP$ a cluster, if it has more than one element itself. Denote the subset of clusters
of $\CP_z$ by $\CC(\CP_z)$.
Now for $C\in\CP_z$ define $z_C$ to be the unique element of $\{z_j:j\in C\}$. 
With this notation,
the proof splits essentially into two major cases of types of collisions which have to be excluded. The one which can be settled most easily is the case of interior collisions, that is, given an initial value $z_0\in\CF_N\Omega$ there exists a point $z\in\partial\CF_N\Omega$, such that the
partition $\CP_z$ has a cluster $C\in\CC(\CP)$ satisfying $z_C\in\Omega$, such that $\phi(z_0,t)_j\to z_C$ as $t\to t^+(z_0)$ if and only if $j\in C$. We denote the set of interior collision points as\[
\partial_{\mathrm{int}}\CF_N\Omega=
\left\{z\in\partial\CF_N\Omega:\ \ex\ C\in\CC(\CP_z)\ \text{such that}\ z_C\in\Omega\right\}.
\]
Note that this does include the case of vortices colliding with the boundary, as long as there are some other vortices, which collide inside $\Omega$ at the same time. 

The second case, which in the following is termed "boundary collisions" is more complicated to settle. In this case the collision point $z\in\partial\CF_N\Omega$ satisfies $z_C\in\partial\Omega$ for each cluster $C\in\CC(\CP)$, and it holds that $\phi(z_0,t)_j\to z_C$ as $t\to t^+(z_0)$ for all $j\in C$, $C\in\CC(\CP_z)$.
We denote the set of boundary collisions with\[
\partial_{\mathrm{bdry}}\CF_N\Omega=\left\{z\in\partial\CF_N\Omega:\ \fa\ C\in\CC(\CP_z): z_C\in\partial\Omega\right\}.
\]
Clearly $\partial\CF_N\Omega$ is the disjoint union of these two sets.

Before we turn to the proof of proposition \ref{nocollisions}, we state some essential lemmata which help us settle the above two cases.

\subsection{Interior collisions}

\begin{lemma}\label{intr_help_fn}
Let $\Gamma$ be $\Delta$--admissible and for any partition $\CP$ of $\{1,\dots,N\}$, $C\in\CC(\CP)$ define \[
J_C:\CF_N\IC\ni z\mapsto \sum\limits_{\substack{i,j\in C\\ i\not=j}}\Gamma_i\Gamma_j\ln|z_i-z_j|\in\IR,
\]\[
J_{\CP}:\CF_N\IC\ni z\mapsto \sum\limits_{C\in\CC(\CP)}J_C(z)\in\IR.        
\]
Further define the constant $C_\Gamma$ by \[
C_\Gamma:=\min\limits_{\substack{\CP\, \mathrm{partition}\\ \mathrm{of}\, \{1,\dots,N\}\\ \CC(\CP)\not=\emptyset}}\ \min\limits_{C\in\CC(\CP)}\left|\sum\limits_{\substack{i,j\in C\\i\not=j}}\Gamma_i\Gamma_j\right|,
\]
which is positive since $\Gamma$ is $\Delta$--admissible. With this notation the inequality
\[
|\nabla J_{\CP_z}(w)|\ge C_\Gamma\max\limits_{C\in\CC(\CP_z)}\left(\sum\limits_{i\in C}|w_i-z_C|^2\right)^{-\frac{1}{2}}.
\]
holds for every $z\in\partial\CF_N\IC$, $w\in\CF_N\IC$.
%
\end{lemma}

\begin{proof2} 
Fix points $z\in\partial\CF_N\IC$, $w\in\CF_N\IC$ and some cluster $C\in\CC(\CP_z)$, put $\widetilde{z_C}:=(z_C,\dots,z_C)\in\IC^N$ and define
\[
j_C:(0,\infty)\ni r\mapsto J_C(\widetilde{z_C}+r(w-\widetilde{z_C}))\in\IR.
\]
Then \[
j_C'(1)=\sum\limits_{\substack{i,j\in C\\i\not=j}}\Gamma_i\Gamma_j\not=0
\]
and letting $c_C:=|j_C'(1)|$ we infer\[
0<c_C=|j_C'(1)|=\left|\left\langle \nabla J_C(w),w-\widetilde{z_C}\right\rangle\right|\le |\nabla J_C(w)|\cdot \left(\sum\limits_{i\in C}|w_i-z_C|^2\right)^{\frac{1}{2}}
\]
for any $w\in\CF_N\IC$.
Together with $\min\limits_{C\in\CC(\CP_z)}c_C\ge C_\Gamma$ and \[
|\nabla J_{\CP_z}(w)|=\left(\sum\limits_{C\in\CC(\CP_z)}\left|\nabla J_C(w)\right|^2\right)^{\frac{1}{2}}\ge \max\limits_{C\in\CC(\CP_z)}|\nabla J_C(w)|
\]
the claim follows. 
\end{proof2}

\begin{lemma}\label{interior_collision}
Let $\Gamma$ be $\Delta$--admissible and $\bar z\in\partial_{\mathrm{int}}\CF_N\Omega$ with corresponding partition $\CP_{\bar z}$, and let $C\in\CC(\CP_{\bar z})$ be an interior collision cluster, that is $\bar z_C\in\Omega$.
There exists $\delta>0$, such that for each $z\in U_\delta(\bar z)\cap\CF_N\Omega$:\[
|\nabla H_\Omega(z)|\ge  \frac{C_\Gamma}{4\pi}\left(\sum\limits_{i\in C}|z_i-\bar z_C|^2\right)^{-\frac{1}{2}}.
\]
\end{lemma}

\begin{proof2}We decompose
$H_\Omega$ as\[
H_\Omega(z)=-\frac{1}{2\pi}J_{\CP_{\bar z}}(z)+K(z),
\]
where
\[K(z)=\sum\limits_{j=1}^N\Gamma_j^2h(z_j)+\sum\limits_{I\in\CC(\CP_{\bar z})}\sum\limits_{\substack{i,j\in I\\ i\not=j}}\Gamma_i\Gamma_jg(z_i,z_j)
+\sum\limits_{\substack{i,j\in \{1,\dots,N\}\\ \ex I,J\in\CP_{\bar z}, i\in I, j\in J\\ I\cap J=\emptyset}}\Gamma_i\Gamma_j G(z_i,z_j).\]
Fixing some interior collision cluster $C$, we have \[
|\nabla H_\Omega(z)|\ge |\nabla H_\Omega(z)|_C=\left|-\frac{1}{2\pi}\nabla J_{\CP_{\bar z}}(z)+\nabla K(z)\right|_C\]\[
\ge \frac{1}{2\pi}|\nabla J_{\CP_{\bar z}}(z)|_C-|\nabla K(z)|_C=\frac{1}{2\pi}|\nabla J_C(z)|-|\nabla K(z)|_C
\]\begin{equation}\label{splitted_intr_grad}\begin{split}
\ge \frac{1}{2\pi}|\nabla J_C(z)|-\sum\limits_{j\in C}\Gamma_j^2|\nabla h(z_j)|-\sum\limits_{\substack{i,j\in C\\ i\not=j}}|\Gamma_i\Gamma_j|\left(|\nabla_{z_i}g(z_i,z_j)|+|\nabla_{z_j}g(z_i,z_j)|\right)\\
-\sum\limits_{i\in C,j\not\in C}|\Gamma_i\Gamma_j| |\nabla_{z_i} G(z_i,z_j)|,\end{split}
\end{equation}
where for $\zeta\in\IC^N$: $|\zeta|_C:=|\pi_C\zeta|$, and $\pi_C:\IC^N\to\{z\in\IC^N:z_j=0 \ \text{if}\ j\not\in C\}$ is the orthogonal projection. For $z\in\IC^N$ define\[
r_C(z):=\min\left\{\min_{i\in C,j\not\in C}|z_i-z_j|,\min_{j\in C}\dist(z_j,\partial\Omega)\right\}.
\]
Since $\delta_0:=r_C(\bar z)>0$ and $r_C$ is clearly continuous, there is $\tilde\delta>0$, such that $r_C(z)\ge\tfrac{\delta_0}{2}$ for every $z\in U_{\tilde\delta}(\bar z)\cap\CF_N\Omega$, which by means of hyothesis \ref{assumptiononG} implies that on $U_{\tilde\delta}(\bar z)\cap\CF_N\Omega$ the last terms of \eqref{splitted_intr_grad} are bounded by some constant $\tilde C\ge 0$.

Now applying lemma \ref{intr_help_fn} yields\[
|\nabla H_\Omega(z)|\ge  \frac{C_\Gamma}{2\pi}\left(\sum\limits_{i\in C}|z_i-\bar z_C|^2\right)^{-\frac{1}{2}}-\tilde C.
\]
Since $\left(\sum\limits_{i\in C}|z_i-\bar z_C|^2\right)^{-\frac{1}{2}}\to\infty$ for $z\to\bar z$, we may choose some $\delta\in(0,\tilde\delta)$, such that for every $z\in U_\delta(\bar z)\cap\CF_N\Omega$:\[
\left(\sum\limits_{i\in C}|z_i-\bar z_C|^2\right)^{-\frac{1}{2}}\ge \frac{4\pi \tilde C}{C_\Gamma},
\]
so that\[
|\nabla H_\Omega(z)|\ge  \frac{C_\Gamma}{4\pi}\left(\sum\limits_{i\in C}|z_i-\bar z_C|^2\right)^{-\frac{1}{2}},
\]
which is what we were to show.
\end{proof2}

\subsection{Boundary collisions}

Now we study the behaviour of $H_\Omega$ near $\partial_{\mathrm{bdry}}\CF_N\Omega$. Let therefore be $\bar z\in \partial_{\mathrm{bdry}}\CF_N\Omega$, that is
$\CP_{\bar z}$ is a partition of $\{1,\dots,N\}$, such that we have distinct points $\bar z_C\in\partial\Omega$ for every cluster $C\in\CC(\CP_{\bar z})$. It may as well be that $\CC(\CP_{\bar z})=\emptyset$. In this case
we have $\bar z\in\CF_N\partial\Omega$.

Settling the case of interior collisions is relatively easy since, away from $\partial\Omega$, the logarithmic singularity of $G$ dominates the interaction between vortices. If two vortices $x,y$ are near to the boundary and to each other, this is no longer true, since then the term $g(x,y)$ cannot be neglected. The next lemma is the key to the
understanding of the interaction taking place between vortices near the boundary.

\begin{lemma}\label{bdry_interaction}
Setting\[
A(x,y):=2\pi\left(\spr{\nabla_xG(x,y),d_x\nu_x}+\spr{\nabla_yG(x,y),d_y\nu_y}\right)
\]
for $x,y\in\Omega_\eps$, we have\[
A(x,y)=\frac{\spr{x-\bar y,\nu_y}^2}{|x-\bar y|^2}-\frac{\spr{x-y,\nu_y}^2}{|x-y|^2}+o(1),
\]
as well as\[
|A(x,y)|\le 1+o(1)
\]
as $x,y\to x_0\in\partial\Omega$. Moreover, if $\Omega$ is strictly convex, we have\[
o(1)\le A(x,y)\le 1+o(1)
\]
\end{lemma}

\begin{proof2}
By hypothesis \ref{assumptiononG} we may write\[
G(x,y)=\frac{1}{2\pi}\ln\frac{|\bar x-y|}{|x-y|}+O(1),
\]
as $x,y\to x_0\in\partial\Omega$, and the approximation holds in the $C^1$--sense, therefore (since $\bar x=2p(x)-x$)\[
2\pi\nabla_x G(x,y)=(2\Dif p(x)-\id)\frac{\bar x-y}{|\bar x-y|^2}-\frac{x-y}{|x-y|^2}+O(1)
\]\[
=\frac{y-\bar x}{|\bar x-y|^2}-\frac{x-y}{|x-y|^2}+\frac{2}{1-\kappa_xd_x}\spr{\frac{\bar x-y}{|\bar x-y|^2},\tau_x}\tau_x+O(1),
\]
which leads to\[
A(x,y)=\spr{\frac{y-\bar x}{|\bar x-y|^2}-\frac{x-y}{|x-y|^2},d_x\nu_x}+\spr{\frac{x-\bar y}{|\bar y-x|^2}-\frac{y-x}{|x-y|^2},d_y\nu_y}+o(1)
\]\[
=\spr{\frac{y-x+2d_x\nu_x}{|\bar x-y|^2},d_x\nu_x}+\spr{\frac{x-y+2d_y\nu_y}{|\bar y-x|^2},d_y\nu_y}-\spr{\frac{x-y}{|x-y|^2},d_x\nu_x-d_y\nu_y}+o(1)
\]\[\begin{split}
=2\left(\frac{d_x^2}{|\bar x-y|^2}+\frac{d_y^2}{|\bar y-x|^2}\right)+\spr{y-x,\frac{d_x\nu_x}{|\bar x-y|^2}}+\spr{x-y,\frac{d_y\nu_y}{|\bar y-x|^2}}
\\-\spr{\frac{x-y}{|x-y|^2},d_x\nu_x-d_y\nu_y}+o(1)
\end{split}
\]\[\begin{split}
=2\left(\frac{d_x^2}{|\bar x-y|^2}+\frac{d_y^2}{|\bar y-x|^2}\right)+\spr{x-y,\frac{d_y\nu_y}{|\bar y-x|^2}-\frac{d_x\nu_x}{|\bar x-y|^2}}-1\\+\spr{\frac{x-y}{|x-y|^2},p(x)-p(y)}+o(1)
\end{split}
\]
Since $p(x)-p(y)=\Dif p(y)(x-y)+o(|x-y|)$ we have, using lemma \ref{Dp}\[
\begin{split}
A(x,y)=2\left(\frac{d_x^2}{|\bar x-y|^2}+\frac{d_y^2}{|\bar y-x|^2}\right)+\spr{x-y,\frac{d_y\nu_y}{|\bar y-x|^2}-\frac{d_x\nu_x}{|\bar x-y|^2}}-1\\+\frac{\spr{x-y,\tau_y}^2}{(1-\kappa_yd_y)|x-y|^2}+o(1)
\end{split}
\]\begin{equation}\label{step1}
\begin{split}
=2\left(\frac{d_x^2}{|\bar x-y|^2}+\frac{d_y^2}{|\bar y-x|^2}\right)-\spr{x-y,\frac{d_x\nu_x-d_y\nu_y}{|x-\bar y|^2}}+\frac{\spr{x-y,\tau_y}^2}{(1-\kappa_yd_y)|x-y|^2}-1\\-\spr{x-y,{d_x\nu_x}}\left(\frac{1}{|\bar x-y|^2}-\frac{1}{|x-\bar y|^2}\right)+o(1).
\end{split}
\end{equation}
We shall now see, that the whole last line of \eqref{step1} is in fact $o(1)$:\[
\left|\spr{x-y,{d_x\nu_x}}\left(\frac{1}{|\bar x-y|^2}-\frac{1}{|x-\bar y|^2}\right)\right|\le d_x|x-y|\cdot\frac{1}{|\bar x-y|^2|x-\bar y|^2}\cdot\left||\bar x-y|^2-|x-\bar y|^2\right|
\]\[
\le \frac{d_x|x-y|}{|\bar x-y|^2|x-\bar y|^2}\cdot c(d_x+d_y)|p(x)-p(y)|^2\le\tilde c\frac{d_x(d_x+d_y)|\bar x-y|^2}{|\bar x-y|^2|x-\bar y|^2}\cdot|x-y|
\]\[
\le \tilde c|x-y|=o(1)
\]
for some constants $c,\tilde c>0$, where we used lemma \ref{distinequalities} repeatedly.
For the sake of a more readable presentation, we continue by estimating the second term of \eqref{step1} separately.\[
\spr{x-y,\frac{d_x\nu_x-d_y\nu_y}{|x-\bar y|^2}}=\frac{|x-y|^2}{|x-\bar y|^2}-\spr{\frac{x-y}{|x-\bar y|^2},p(x)-p(y)}
\]\[
=\frac{|x-y|^2}{|x-\bar y|^2}-\spr{\frac{x-y}{|x-\bar y|^2},\frac{1}{1-\kappa_yd_y}\spr{x-y,\tau_y}\tau_y+o(|x-y|)}
\]\[=\frac{|x-y|^2}{|x-\bar y|^2}-\frac{\spr{x-y,\tau_y}^2}{(1-\kappa_yd_y)|x-\bar y|^2}+o(1)
\]
Therefore\[
A(x,y)=2\!\left(\frac{d_x^2}{|\bar x-y|^2}+\frac{d_y^2}{|\bar y-x|^2}\!\right)-\frac{|x-y|^2}{|x-\bar y|^2}+\frac{\spr{x-y,\tau_y}^2}{1-\kappa_yd_y}\!\left(\frac{1}{|x-\bar y|^2}+\frac{1}{|x-y|^2}\!\right)-1+o(1)
\]\[\begin{split}
=2\left(\frac{d_x^2}{|\bar x-y|^2}+\frac{d_y^2}{|\bar y-x|^2}\right)-\frac{|x-y|^2}{|x-\bar y|^2}+\spr{x-y,\tau_y}^2\left(\frac{1}{|x-\bar y|^2}+\frac{1}{|x-y|^2}\right)-1\\
+\spr{x-y,\tau_y}^2\left(\frac{1}{|x-\bar y|^2}+\frac{1}{|x-y|^2}\right)\left(\frac{1}{1-\kappa_yd_y}-1\right)
+o(1).
\end{split}
\]
Since $|\bar y-x|\ge \hat c|x-y|$ for some $\hat c>0$ and $d_y=o(1)$, whereas $\kappa_y$ is bounded, the last line of the preceding formula is again $o(1)$. Since $|x-y|^2=\spr{x-y,\tau_y}^2+\spr{x-y,\nu_y}^2$ we may then rewrite $A(x,y)$ as\[
A(x,y)=2\left(\!\frac{d_x^2}{|\bar x-y|^2}+\frac{d_y^2}{|\bar y-x|^2}\!\right)-\frac{|x-y|^2}{|x-\bar y|^2}+\spr{x-y,\tau_y}^2\!\left(\frac{1}{|x-\bar y|^2}+\frac{1}{|x-y|^2}\right)-1+o(1)
\]\[
=2\left(\frac{d_x^2}{|\bar x-y|^2}+\frac{d_y^2}{|\bar y-x|^2}\right)-\spr{x-y,\nu_y}^2\left(\frac{1}{|x-\bar y|^2}+\frac{1}{|x-y|^2}\right)+o(1)
\]\[\begin{split}
=2\frac{d_x^2+d_y^2}{|x-\bar y|^2}-\spr{x-y,\nu_y}^2\left(\frac{1}{|x-\bar y|^2}+\frac{1}{|x-y|^2}\right)
\\ -2d_x^2\left(\frac{1}{|x-\bar y|^2}-\frac{1}{|\bar x-y|^2}\right)
+o(1).
\end{split}
\]
Again, the last line is $o(1)$, since\[
\left|d_x^2\left(\!\frac{1}{|x-\bar y|^2}-\frac{1}{|\bar x-y|^2}\!\right)\!\right|=\frac{d_x^2\left||\bar x-y|^2\!-|x-\bar y|^2\right|}{|\bar x-y|^2|x-\bar y|^2}
\le c\frac{d_x^2|p(x)-p(y)|^2}{|\bar x-y|^2|x-\bar y|^2}\cdot(d_x+d_y)=o(1),
\]
which is also obtained using lemma \ref{distinequalities}. In the following we abbreviate $\alpha=\spr{x-y,\nu_y}$, $\beta=\spr{x-y,\tau_y}$, hence\[
x-y=\alpha\nu_y+\beta\tau_y,
\]\[
x-\bar y=x-y+2d_y\nu_y=(\alpha+2d_y)\nu_y+\beta\tau_y,
\]\[
d_x^2=d_y^2+\Dif d^2(y)(x-y)+\tfrac{1}{2}\Dif^2d^2(y)[x-y,x-y]+o(|x-y|^2)\]\[
=d_y^2+2d_y\spr{x-y,\nu_y}+\spr{x-y,(\id-\Dif p(y))(x-y)}+o(|x-y|^2)\]
\[=d_y^2+2\alpha d_y+\alpha^2+\beta^2-\frac{1}{1-\kappa_yd_y}\beta^2+o(|x-y|^2)
\]\[
=d_y^2+2\alpha d_y+\alpha^2+\beta^2\left(1-\frac{1}{1-\kappa_yd_y}\right)+o(|x-y|^2)
\]\[
=d_y^2+2\alpha d_y+\alpha^2+o(|x-y|^2),
\]
which implies\[
A(x,y)=\frac{4d_y^2+4\alpha d_y+2\alpha^2-\alpha^2}{(\alpha+2d_y)^2+\beta^2}-\frac{\alpha^2}{\alpha^2+\beta^2}+o(1)
=\frac{(\alpha+2d_y)^2}{(\alpha+2d_y)^2+\beta^2}-\frac{\alpha^2}{\alpha^2+\beta^2}+o(1),
\]
which is precisely our first claim.
The second claim also follows easily, since $\tfrac{\xi^2}{\xi^2+\beta^2}\in[0,1]$ for every $\xi\in\IR$.

We will now show that the above is in fact nonnegative up to an error of $o(1)$ if $\Omega$ is strictly convex. 

First observe that $A(x,y)=o(1)$ for $\beta=0$. On the other hand, setting $f(t):=\frac{t^2}{t^2+\beta^2}$ for $\beta\not=0$ and $t\in\IR$ we may apply the mean value theorem to get\[
A(x,y)=f(\alpha+2d_y)-f(\alpha)+o(1)=2d_yf'(\xi)+o(1)=\frac{4d_y\beta^2\xi}{(\xi^2+\beta^2)^2}+o(1)
\]
for some $\xi\in[\alpha,\alpha+2d_y]$. 

Now define $\alpha=\spr{x-y,\nu_y}$ and $\alpha':=\spr{y-x,\nu_x}$. In the sequel, we will show that for any $x,y\in\Omega_\eps$ one of both scalar products is nonnegative if $|x-y|$ is small enough. 

Assume on the contrary that there are sequences $(x^n), (y^n)\subset\Omega_\eps$ with $|x^n-y^n|\to 0$ as $n\to\infty$ such that $\alpha_n=\spr{x^n-y^n,\nu_{y^n}}<0$ and $\alpha_n'=\spr{y^n-x^n,\nu_{x^n}}<0$ for every $n\in\IN$. Then we have\[
0>\alpha_n+\alpha'_n=\spr{x^n-y^n,\nu_{y^n}-\nu_{x^n}}=\spr{x^n-y^n,-\Dif\nu_{y^n}(x^n-y^n)+o(|x^n-y^n|)}\]
\[=\spr{x^n-y^n,\frac{\kappa_{y^n}}{1-\kappa_{y^n}d_{y^n}}\spr{x^n-y^n,\tau_{y^n}}\tau_{y^n}+o(|x^n-y^n|)}
\]\[
=\frac{\kappa_{y^n}}{1-\kappa_{y^n}d_{y^n}}\beta_n^2+o(|x^n-y^n|^2),
\]
hence $\beta_n^2:=\spr{x^n-y^n,\tau_{y^n}}^2=o(|x^n-y^n|^2)$, since $\kappa_y\ge\tilde\eps$ for all $y\in\Omega_\eps$ and some $\tilde\eps>0$ if $\Omega$ is strictly convex. Since $\alpha_n^2+\beta_n^2=|x^n-y^n|^2$ we infer 
\[
\frac{\alpha_n^2}{|x^n-y^n|^2}=1-\frac{\beta_n^2}{|x^n-y^n|^2}\to 1
\]
as $n\to \infty$, which implies
\[
\frac{\alpha_n}{|x^n-y^n|}\to -1
\]
as $n\to\infty$, since we have assumed $\alpha_n<0$. It follows that
\[
-1\leftarrow\frac{\alpha_n}{|x^n-y^n|}>\frac{\alpha_n+\alpha'_n}{|x^n-y^n|}=\frac{\kappa_{y^n}}{1-\kappa_{y^n}d_{y^n}}\frac{\beta_n^2}{|x^n-y^n|}+o(|x^n-y^n|)\to 0
\]
as $n\to\infty$, which is the desired contradiction.

Hence for every $x,y\in\Omega_\eps$ sufficiently close to each other, one of the two scalar products $\alpha$, $\alpha'$ is nonnegative, and since $A(x,y)$ is symmetric in $x,y$ by definition, we might interchange the roles of $x$ and $y$ to assume $\alpha\ge0$, which in turn implies $\xi\ge 0$ and we are done.
\end{proof2}

\begin{lemma}\label{nocollisionwithbdry}
Let $\Gamma$ be $\partial$--admissible, $\bar z\in\partial_{\text{bdry}}\CF_N\Omega$ and let $C\in\CP_{\bar z}$ satisfy $\bar z_C\in\partial\Omega$. There is $\delta>0$ such that for every $z\in U_\delta(\bar z)\cap\CF_N\Omega$\[
|\nabla H_\Omega(z)|\ge \frac{\eps_C}{2\pi}\left(\sum_{j\in C}d(z_j)^2\right)^{-\frac{1}{2}},
\]
where the constant $\eps_C>0$ is given by\[\eps_C:=\frac{1}{2}\left(\sum_{i\in C}\Gamma_i^2-\sum_{\substack{i,j\in C\\ i\not=j}}|\Gamma_i\Gamma_j|\right),\]
if $\Omega$ is not strictly convex, and by\[
\eps_C:=\frac{1}{2}\left(\sum_{i\in C}\Gamma_i^2-\sum_{\substack{i,j\in C\\ \Gamma_i\Gamma_j<0}}|\Gamma_i\Gamma_j|\right)
\]
instead, if it is strictly convex.
\end{lemma}
\begin{proof2}
Note that in each case $\eps_C>0$ by the condition of $\partial$--admissibility. There is $\tilde\delta>0$, such that $z_j\in\Omega_\eps$ for any $j\in C$, $z\in U_{\tilde\delta}(\bar z)\cap\CF_N\Omega$. We thus may consider the function\[
\Phi_C:U_{\tilde\delta}(\bar z)\cap\CF_N\Omega\ni z\mapsto\pi\sum_{j\in C}d(z_j)^2\in[0,\infty)
\]
and simply compute\[
\spr{\nabla H_\Omega(z),\nabla\Phi_C(z)}=2\pi\sum_{j\in C}\spr{\nabla_{z_j} H_\Omega(z),d(z_j)\nu(z_j)}
\]\[
=2\pi\sum_{j\in C}\spr{\Gamma_j^2\nabla h(z_j)+2\sum_{i\in C\setminus\{j\}}\Gamma_i\Gamma_j\nabla_{z_j}G(z_j,z_i)+O(1),d(z_j)\nu(z_j)}
\]\[
=2\pi\sum_{j\in C}\spr{\Gamma_j^2\nabla h(z_j),d(z_j)\nu(z_j)}+4\pi\sum_{\substack{i,j\in C\\ i\not=j}}\Gamma_i\Gamma_j\spr{\nabla_{z_j}G(z_j,z_i),d(z_j)\nu(z_j)}+o(1)
\]\[
=\sum_{j\in C}\Gamma_j^2\spr{\frac{\nu(z_j)}{d(z_j)},d(z_j)\nu(z_j)}+2\sum_{\substack{i,j\in C\\ i<j}}\Gamma_i\Gamma_jA(z_i,z_j)+o(1)
\]
as $z\to\bar z$. Here we have used the fact that $\nabla_{z_j} G(z_j,z_i)=O(1)$ for $j\in C$, $i\not\in C$ as $z\to\bar z$ by hypothesis \ref{assumptiononG}.

In case $\Omega$ is not strictly convex, we may estimate this by
\[
\ge\sum_{j\in C}\Gamma_j^2-2\sum_{\substack{i,j\in C\\ i<j}}|\Gamma_i\Gamma_j|+o(1)=\sum_{j\in C}\Gamma_j^2-\sum_{\substack{i,j\in C\\ i\not=j}}|\Gamma_i\Gamma_j|+o(1)=2\eps_C+o(1)
\] 
as $z\to \bar z$ by  use of lemma \ref{bdry_interaction}.

If, on the other hand, $\Omega$ is strictly convex, we may again use lemma \ref{bdry_interaction} and similarly conclude\[
\spr{\nabla H_\Omega(z),\nabla\Phi_C(z)}\ge \sum_{j\in C}\Gamma_j^2-2\sum_{\substack{i,j\in C\\ i<j,\, \Gamma_i\Gamma_j<0}}|\Gamma_i\Gamma_j|+o(1)\]\[=\sum_{j\in C}\Gamma_j^2-\sum_{\substack{i,j\in C\\ \Gamma_i\Gamma_j<0}}|\Gamma_i\Gamma_j|+o(1)=2\eps_C+o(1)
\]
as $z\to\bar z$. In any case we obtain that there is $\delta\in(0,\tilde\delta)$ such that\[
\spr{\nabla H_\Omega(z),\nabla\Phi_C(z)}\ge\eps_C
\]
for every $z\in U_\delta(\bar z)\cap\CF_N\Omega$.

On the other hand we have\[
\spr{\nabla H_\Omega(z),\nabla\Phi_C(z)}\le |\nabla H_\Omega(z)|\cdot|\nabla\Phi_C(z)|=2\pi\cdot |\nabla H_\Omega(z)|\cdot\left(\sum_{j\in C}d(z_j)^2\right)^{\frac{1}{2}}
\]
by simply applying the Cauchy--Schwarz--inequality, hence we obtain\[
|\nabla H_\Omega(z)|\ge\frac{\eps_C}{2\pi}\left(\sum_{j\in C}d(z_j)^2\right)^{-\frac{1}{2}}
\]
for every $z\in U_\delta(\bar z)\cap\CF_N\Omega$ and we are done.
\end{proof2}

\subsection{Proof of proposition \ref{nocollisions}}
Equipped with these estimates we now turn to the proof of proposition \ref{nocollisions}, which is comprised of the next few lemmata.

\begin{lemma}\label{M_delta}
There is $\delta>0$ such that $|\nabla H_\Omega(z)|>1$ for every $z\in\CM_\delta$.
\end{lemma}
\begin{proof2}
Assume on the contrary that there are sequences $\delta_n\to 0$, $\delta_n>0$ and $z^n\in\CM_{\delta_n}$ such that $|\nabla H_\Omega(z^n)|\le 1$ for all $n\in\IN$. Then, upon choosing a convergent subsequence, we may assume $z^n\in U_{\delta_n}(\bar z)$ for some $\bar z\in\partial\CF_N\Omega$ and every $n\in\IN$. Now if $\bar z\in\partial_{\mathrm{int}}\CF_N\Omega$ there are $n_0\in\IN$ and $C\in\CC(\CP_{\bar z})$ satisfying $\bar z_C\in\Omega$ such that \[
|\nabla H_\Omega(z^n)|\ge \frac{C_\Gamma}{4\pi}\left(\sum\limits_{i\in C}|z_i^n-\bar z_C|^2\right)^{-\frac{1}{2}}
\ge \frac{C_\Gamma}{4\pi\sqrt{|C|}\delta_n}
\] for every $n\ge n_0$ by lemma \ref{interior_collision}. Similarly, if $\bar z\in\partial_{\mathrm{bdry}}\CF_N\Omega$ there are $n_0\in\IN$ and $C\in\CP_{\bar z}$ satisfying $\bar z_C\in\partial\Omega$ such that\[
|\nabla H_\Omega(z^n)|\ge \frac{\eps_C}{2\pi}\left(\sum_{j\in C}d(z_j^n)^2\right)^{-\frac{1}{2}}\ge\frac{\eps_C}{2\pi\sqrt{|C|}\delta_n}
\]
for every $n\ge n_0$ by lemma \ref{nocollisionwithbdry}. This, however, contradicts the fact that $|\nabla H_\Omega(z^n)|\le 1$ for all $n\in\IN$ and the proof is done.
\end{proof2}

Since $\overline{\CM_\delta}$ is a neighbourhood of $\partial\CF_N\Omega$ in $\overline{\Omega}^N$, this in particular shows that $H_\Omega$ satisfies the Palais--Smale--condition.

\begin{lemma}\label{t+infty}
Let $z\in\CF_N\Omega$ satisfy\[
\lim\limits_{t\to t^+(z)}H_\Omega(\phi(z,t))=c_0<\infty.
\]
Then $t^+(z)=\infty$.
\end{lemma}

\begin{proof2}
In general we have for $s,t\in[0,t^+(z))$, $s<t$ and $z^t:=\phi(z,t)$\[
|z^t-z^s|\le\int\limits_s^t|\nabla H_\Omega(z^\tau)|\dif\tau\le\sqrt{t-s}\sqrt{\int_s^t|\nabla H_\Omega(z^\tau)|^2\dif\tau}
=\sqrt{t-s}\sqrt{H_\Omega(z^t)-H_\Omega(z^s)}
\]\begin{equation}\label{cauchy}
\le\sqrt{|t-s|}\sqrt{c_0-H_\Omega(z^s)}  
\end{equation}
Now if $t^+(z)<\infty$ we may take the limit $t\to t^+(z)$ on the right hand side of \eqref{cauchy} to obtain that for every $\eps>0$ there is $t_0\in[0,t^+(z))$ and
any $s,t\in[t_0,t^+(z))$, $s<t$: $|z^t-z^s|<\eps$, hence $z^t\to \bar z$ as $t\to t^+(z)$ for some $\bar z\in\overline{\CF_N\Omega}=\overline{\Omega}^N$, since $\overline{\Omega}^N$ is compact. 

Let $\delta>0$ be such that the consequences of lemmata \ref{interior_collision} and \ref{nocollisionwithbdry} hold. 

If $\bar{z}\in\partial_{\mathrm{bdry}}\CF_N\Omega$ we find $C\in\CP_{\bar{z}}$ and $t_0\in[0,t^+(z))$ such that for every $t\in[t_0,t^+(z))$\[
\frac{\dif}{\dif t}\Phi_C(z^t)=\spr{\nabla H_\Omega(z^t),\nabla\Phi_C(z^t)}\ge\eps_C>0,
\]
by application of lemma \ref{nocollisionwithbdry} which is a contradiction.

If, on the other hand, $\bar z\in\partial_\mathrm{int}\CF_N\Omega$, we have $C\in\CC(\CP_{\bar{z}})$ as well as $t_0\in[0,t^+(z))$, such that for all $t\in[t_0,t^+(z))$\[
|\nabla H_\Omega(z^t)|\ge \frac{C_\Gamma}{4\pi}\left(\sum\limits_{i\in C}|z_i^t-\bar z_C|^2\right)^{-\frac{1}{2}}
\]
by lemma \ref{interior_collision}. We thus may compute for $s\in[t_0,t^+(z))$,
$t\in(s,t^+(z))$\[
H_\Omega(z^t)-H_\Omega(z^s)\ge \int\limits_s^t \frac{C_\Gamma}{4\pi}|\nabla H_\Omega(z^\tau)|\left(\sum\limits_{i\in C}|z_i^\tau-\bar z_C|^2\right)^{-\frac{1}{2}}\dif\tau\]\[
=\frac{C_\Gamma}{4\pi}\int\limits_s^t|\dot z^\tau|\left(\sum\limits_{i\in C}|z_i^\tau-\bar z_C|^2\right)^{-\frac{1}{2}}\dif\tau
\ge \frac{C_\Gamma}{4\pi}\int\limits_s^t|\pi_C{\dot z^\tau}|\left(\sum\limits_{i\in C}|z_i^\tau-\bar z_C|^2\right)^{-\frac{1}{2}}\dif\tau
\]
\[\ge \frac{C_\Gamma}{4\pi}\int\limits_s^t\left|\left\langle\pi_C{\dot z^\tau},\frac{\pi_C(z^\tau-\bar z)}{\left|\pi_C({z^\tau}-{\bar z})\right|}\right\rangle\right|\left(\sum\limits_{i\in C}|z_i^\tau-\bar z_C|^2\right)^{-\frac{1}{2}}\dif\tau\]\[
=\frac{C_\Gamma}{4\pi}\int\limits_s^t\frac{\left|\frac{\dif}{\dif\tau}\big(\sum_{i\in C}|z_i^\tau-\bar z_C|^2\big)^{\frac{1}{2}}\right|}{\big(\sum_{i\in C}|z_i^\tau-\bar z_C|^2\big)^{\frac{1}{2}}}\ \dif\tau \ge -\frac{C_\Gamma}{4\pi}\int\limits_s^t\frac{\frac{\dif}{\dif\tau}\big(\sum_{i\in C}|z_i^\tau\!-\bar z_C|^2\big)^{\frac{1}{2}}}
{\big(\sum_{i\in C}|z_i^\tau-\bar z_C|^2\big)^{\frac{1}{2}}}\ \dif\tau
\]\[
=\frac{C_\Gamma}{8\pi}\ln\frac{\sum_{i\in C}|z_i^s-\bar z_C|^2}
{\sum_{i\in C}|z_i^t-\bar z_C|^2}\to\infty \quad (t\to t^+(z)),
\]
contrary to our assumption. It follows that $t^+(z)=\infty$, which is what we were to show.
\end{proof2}

The next lemma finishes the proof of proposition \ref{nocollisions}.

\begin{lemma}
If there is $z\in\CF_N\Omega$ satisfying\[
\lim\limits_{t\to t^+(z)}H_\Omega(\phi(z,t))=c_0<\infty,
\]
there is a sequence $(z^n)\subset\CF_N\Omega$ and a point $z^*\in\CF_N\Omega$ such that $z^n\to z^*$ and $\nabla H_\Omega(z^n)\to 0$ as $n\to\infty$.
\end{lemma}

\begin{proof2}
Since $t^+(z)=\infty$ by lemma \ref{t+infty}, consider a sequence $(t_n)\subset[0,\infty)$, $t_n\to\infty$, such that $z^{t_n}\to z^*$ as $n\to\infty$. 

Let us assume at first that there are $n_0\in\IN$, $\delta>0$, such that $\delta\le t_{n+1}-t_n$ and $|\nabla H_\Omega(z^t)|\ge\frac{1}{\sqrt{n}}$ for all $n\ge n_0$ and for all $t\in[t_n,t_n+\delta]$. Then we have\[
H_\Omega(z^{t_n})\ge H_\Omega(z^{t_{n_0}})+\sum\limits_{j=n_0}^n\int\limits_{t_j}^{t_j+\delta}|\nabla H_\Omega(z^s)|^2\dif s=H_\Omega(z^{t_{n_0}})+
\sum\limits_{j=n_0}^n\frac{\delta}{j}\to\infty \quad (n\to\infty),
\]
contrary to our main assumption. 

Thus there exists a sequence $\delta_n\to 0$, $\delta_n>0$ for all $n\in\IN$, such that with $s_n:=t_n+\delta_n$: $|\nabla H_\Omega(z^{s_n})|<\frac{1}{\sqrt{n}}$.
It follows that $z^{s_n}\to z^*$ as $n\to\infty$  because of \eqref{cauchy}, since $|s_n-t_n|=\delta_n\to 0$. Abbreviating $z^{s_n}$ by $z^n$ we have $\nabla H_\Omega(z^n)\to 0$ as well as $z^n\to z^*$ as $n\to\infty$. Lemma \ref{M_delta} now implies $z^*\in\CF_N\Omega$ and the proof is finished.
\end{proof2}

This also finishes the proof of proposition \ref{nocollisions}. All that is left to do for proving our main results concerning asymmetric domains is to provide a sort of linking argument for $H_\Omega$, that is finding a point $z\in\CF_N\Omega$ such that $H_\Omega(z^t)$ has a finite limit for $t\to\infty$. This is done by applying lemma
\ref{deformation} within the next section.

\section{Linking phenomena for \texorpdfstring{$H_\Omega$}{the functional}}\label{linking_arguments}

\subsection{The simply connected case}
This subsection is concerned with the proof of theorem \ref{general}. 

Let $\Gamma$ be $\partial$-- and $\Delta$--admissible, and let $\Gamma$ be $\CL$--admissible with corresponding
permutation $\sigma$. Reordering the vortices, we may without loss of generality assume that $\sigma=\id$ and hence abbreviate $\CL_N\Omega:=\CL_N^{\id}\Omega$.

The theorem is trivial for $N\le 2$, for then $H_\Omega(z)\to-\infty$ for $z\to\partial\Omega$ and consequently $H_\Omega$ assumes a local maximum in $\CF_N\Omega$, since $h(z_j)\to-\infty$ for $z_j\to\partial\Omega$, $j\in\{1,2\}$ as well as $\Gamma_1\Gamma_2G(z_1,z_2)\to-\infty$
if $|z_1-z_2|\to 0$, since $\Gamma_1\Gamma_2<0$.

In what comes we thus consider the case $N\ge 3$ and begin to construct an explicit linking for $H_\Omega$.

Without loss of generality we may assume $0\in\Omega$. Choose $\rho>0$ such that $\overline{B_{N\rho}(0)}\subset\Omega$. 
Define\[
\gamma_0:\IT^{N-2}:=(S^1)^{N-2}\to\CF_N\Omega 
\]\[
\gamma_{0,1}(\zeta_1,\dots,\zeta_{N-2}):= 0, \qquad \gamma_{0,N}(\zeta_1,\dots,\zeta_{N-2}):= (N-1)\rho, 
\]\[ 
\gamma_{0,j}(\zeta_1,\dots,\zeta_{N-2})=(j-1)\rho\zeta_{j-1},
\]
for $j\in\{2,\dots,N-1\}$, where $\gamma_{0,j}$ denotes the $j$-th component of $\gamma_0$. Setting\[
\Gamma_0:=\left\{\gamma\in C^0(\IT^{N-2},\CF_N\Omega): \gamma\simeq\gamma_0\right\},
\]
we have the following
\begin{lemma}
For every $\gamma\in\Gamma_0$: $\gamma(\IT^{N-2})\cap\CL_N\Omega\not=\emptyset$.
\end{lemma}

\begin{figure}[ht]
\begin{center}
\includegraphics[scale=1]{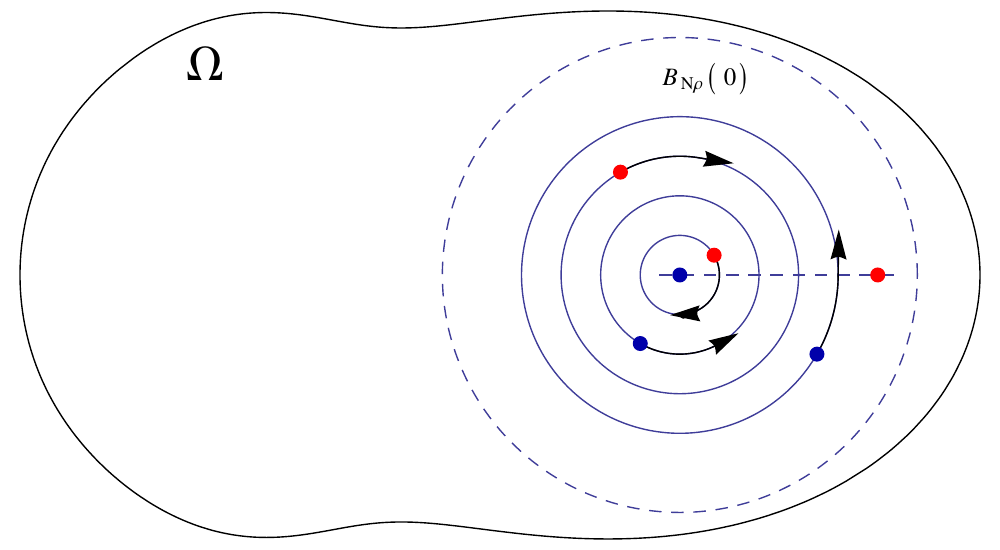}
\end{center}
\caption{The action of the map $\gamma_0$: Vortices are rotated in concentric circles around the point $0\in\Omega$.}
\end{figure}

\begin{proof2}
Let $\tilde H:\IT^{N-2}\times[0,1]\to\CF_N\Omega$ be a deformation from $\gamma_0$ to $\gamma$.
For $t\in[0,1]$ define\[
h_t:\IT^{N-2}\times[0,1]^{N-2}\to \IC^{N-2}
\]
by setting
\[
h_{t,j}(\zeta,s):=s_j\left(\tilde H_j(\zeta,t)-\tilde H_{j+1}(\zeta,t)\right)+(1-s_j)\left(\tilde H_{j+2}(\zeta,t)-\tilde H_{j+1}(\zeta,t)\right)
\]
for $j\in\{1,\dots,N-2\}$.
Obviously $\gamma(\zeta)\in\CL_N\Omega$ if and only if $h_1(\zeta,s)=0$ for some $s\in[0,1]^{N-2}$. 
Furthermore $h_t(\zeta,s)\not=0$ for all $s\in\partial\left([0,1]^{N-2}\right)$, $t\in[0,1]$,
since $\tilde H(\zeta,t)\in\CF_N\Omega$, so the map\begin{equation}\label{pairs}
\!\!\!\begin{array}{c}
\tilde g:\!\left(\IT^{N-2}\!\times[0,1]^{N-2}\!\times[0,1],\IT^{N-2}\!\times\partial\left([0,1]^{N-2}\right)\times[0,1]\right)
\to\left(\IC^{N-2},\!\IC^{N-2}\!\setminus\!\{0\}\right) \\
\tilde g(\zeta,s,t):= h_t(\zeta,s)
\end{array}
\end{equation}
is well--defined and continuous.

Using the K\"{u}nneth--formula for the pair 
\[(X,A):=\left(\IT^{N-2}\times[0,1]^{N-2},\IT^{N-2}\times\partial\left([0,1]^{N-2}\right)\right)\cong\left(\IT^{N-2}\times D^{N-2},\IT^{N-2}\times S^{N-3}\right)\]
we easily get that 
\[H_{2N-4}\left(\IT^{N-2}\times[0,1]^{N-2},\IT^{N-2}\times\partial\left([0,1]^{N-2}\right)\right)\cong\IZ,\]
where we are using singular homology with coefficients in $\IZ$.
Since $\tilde{g}$ is a homotopy of pairs by \eqref{pairs}, $h_t$ induces a homomorphism in homology\[
h_*:H_{2N-4}\left(\IT^{N-2}\times[0,1]^{N-2},\IT^{N-2}\times\partial\left([0,1]^{N-2}\right)\right)
\to H_{2N-4}\left(\IC^{N-2},\IC^{N-2}\setminus\{0\}\right),
\]
which is independent of $t\in[0,1]$. We claim that the degree $h_*(1)\in\IZ$ of $h_0$ is nonzero.

Observe that $h_0$ has a unique zero at $p:=(\zeta_0,s_0):=\left(1,\dots,1,\tfrac{1}{2},\dots,\tfrac{1}{2}\right)$. 
Abbreviating $(Y,B):=\left(\IC^{N-2},\IC^{N-2}\setminus\{0\}\right)$ we thus have the following commutative diagram:
\begin{diagram}
 (X,A) & \rTo^{h_0} & (Y,B)\\
 \dInto^{i} & \ruTo^{u}  & \uTo_{w=h_0|_{(U_\eps(p),U_\eps(p)\setminus\{p\})}} \\
 (X,X\setminus\{p\}) & \lInto^{j} & (U_\eps(p),U_\eps(p)\setminus\{p\}),
\end{diagram}
where $u$ in the middle is given by $u:(X,X\setminus\{p\})\to (Y,B),x\mapsto h_0(x)$.
Taking the $(2N-4)$-th homology, we notice that the restriction homomorphism $i_*$ is an isomorphism since $A$ is the boundary of the $\partial$--manifold $X$, where $X\setminus A$ is orientable and $i_*$ maps a generator $\{X\}$ of $H_{2N-4}(X,A)$, which is a fundamental class corresponding to a global orientation of $X$ to a local orientation of $X$, i.e. a generator of the local homology group $H_{2N-4}(X,X\setminus\{p\})$ of $X$. See, for example \cite{tomdieck}, chapter V, theorem 13.1 for further details and rigorous proofs.
Since $j_*$ is an excision isomorphism, we have \[
h_*=u_*i_*=w_*j_*^{-1}i_*,                                                 
\]
so we are done if the map $w_*$ to the right is an isomorphism. But this is surely the case, as 
$w_*(1)\in\IZ$ for small $\eps>0$ is the local degree of the differentiable map $h_0$ at $p$ and is nonzero,
which can be easily computed as follows:
We have (regarding $\IC^{N-2}$ as $\IR^{2N-4}$)
\[
\frac{\partial}{\partial\zeta_{j-1}}h_{0,j}(p)=\frac{1}{2}(j-1)\rho\begin{pmatrix}0\\ 1\end{pmatrix}, \ \
\frac{\partial}{\partial\zeta_{j}}h_{0,j}(p)=-j\rho\begin{pmatrix}0\\ 1\end{pmatrix}, \ \
\frac{\partial}{\partial\zeta_{j+1}}h_{0,j}(p)=\frac{1}{2}(j+1)\rho\begin{pmatrix}0\\ 1\end{pmatrix},
\]\[
\frac{\partial}{\partial s_{j}}h_{0,j}(p)=-2\rho\begin{pmatrix}1\\ 0\end{pmatrix},
\]
whereas all the other partial derivatives vanish.
Reordering the Jacobian of $h_0$ at $p$, such that the first $N-2$ rows correspond to the imaginary parts of the $h_{0,j}$, we get that\[
\Dif h_0(p)=\begin{pmatrix}A_N & 0\\ 0 & B_N\end{pmatrix},
\]
where $B_N=\diag(-2\rho)\in\IR^{(N-2)\times(N-2)}$. Developing the last $N-2$ columns of\\ $\det\Dif h_0(p)$, and using multilinearity to get rid of the $\rho$-factors, we get
\[
|\det\Dif h_0(p)|=2^{N-2}\rho^{2N-4}\operatorname{abs}\begin{vmatrix}
-1 & 1 & 0 & \cdots & \cdots &\cdots & 0\\
\tfrac{1}{2} & -2 & \tfrac{3}{2} & 0 & & & \vdots \\
0 & 1 & -3 & 2 & 0 & & \vdots \\
\vdots & \ddots & \ddots & \ddots & \ddots & \ddots & \vdots\\
\vdots &   & 0 & \tfrac{N-5}{2} & \scriptstyle{-(N-4)} & \tfrac{N-3}{2} & 0\\
\vdots &   &   &   0   &  \tfrac{N-4}{2} & \scriptstyle{-(N-3)} & \tfrac{N-2}{2} \\
0 & \cdots & \cdots & \cdots & 0 & \tfrac{N-3}{2} & \scriptstyle{-(N-2)}
\end{vmatrix}
\]
By looking at the first $N-3$ rows of the above matrix, we see that if $v\in\IR^{N-2}$ is in the kernel, $v$ must have all components equal to some $v_0\in\IR$. The last row 
then implies that $v_0=0$, hence $\det\Dif h_0(p)\not=0$ and we are done.
\end{proof2}

By applying lemma \ref{deformation} to the homotopy class $\Gamma_0$, we get that there is $\zeta\in\IT^{N-2}$ such that $H_\Omega$ remains bounded along the trajectory of
$\gamma_0(\zeta)$. Since $H_\Omega$ satisfies the Palais--Smale--condition, the proof of theorem \ref{general} is finished.

\newpage

\subsection{\texorpdfstring{The case $\pi_1(\Omega)\not=0$}{The multiple--connected case}}
This subsection is devoted to the proof of theorem \ref{notsimplyconnected}. Hence let the paramter $\Gamma$ be $\partial$-- and $\Delta$--admissible.

Let $\pi_1(\Omega)\not=0$, that is $k_0=\rank\pi_1(\Omega)\ge 1$ and select a bounded component $\Omega_1$ of $\IC\setminus\bar\Omega$. Without
loss of generality we may assume $0\in\Omega_1$.
Defining\[
\CS:=\left\{z\in\CF_N\Omega:\ \fa j\in\{1,\dots,N\}: 
\frac{z_j}{|z_j|}=e^{\frac{2\pi ij}{N}}\right\},
\]
we have the following
\begin{lemma}\label{starbound}
$H_\Omega\big|_\CS$ is bounded from above.
\end{lemma}
\begin{proof2}
Hypothesis \ref{assumptiononOmega} implies that $\Omega$ satisfies an exterior ball condition, hence there is $\rho>0$ such that $|z_i-z_j|>\rho$ for $i,j\in\{1,\dots,N\}$, $i\not=j$. Using hypothesis \ref{assumptiononG}, the assertion is immediate.
\end{proof2}

Since $\bar\Omega$ is a $\partial$--manifold with $\partial\Omega_1$ a component
of $\partial\Omega$ there exists a collar of $\partial\Omega_1\cong S^1$, that is, an open neighborhood $U$ of $\partial\Omega_1$ in $\bar\Omega$ and a homeomorphism
\[\tilde h:S^1\times[0,1)\to U\]
satisfying $h(S^1\times\{0\})=\partial\Omega_1$. Setting \[
h_j:=\tilde h|_{S^1\times\left\{\frac{j}{N+1}\right\}}:S^1\to\Omega
\]
for $j\in\{1,\dots,N\}$, the $h_j$ are Jordan curves with disjoint images enclosing $\Omega_1$. Thus
\[
\gamma_0:=h_1\times\dots\times h_N:\IT^N\to\CF_N\Omega
\]
is well defined, and setting\[
\Gamma_0:=\{\gamma:\IT^N\to\CF_N\Omega:\gamma\simeq\gamma_0\},
\]
we have the following
\begin{lemma}\label{intersection_pi1}
For all $\gamma\in\Gamma_0$\[
\gamma\left(\IT^N\right)\cap\CS\not=\emptyset.
\]
\end{lemma}
\begin{figure}[ht]
\begin{center}
\includegraphics[scale=1]{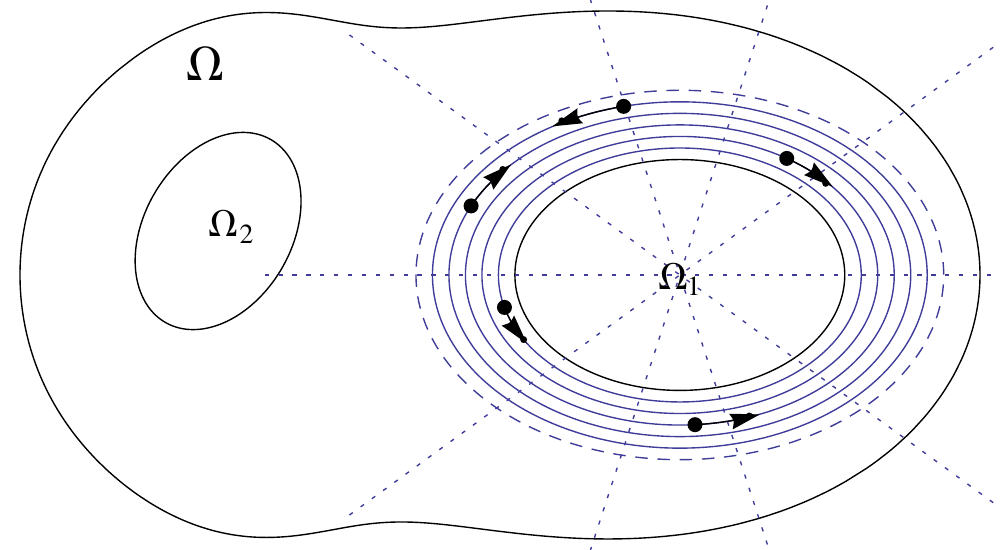}
\end{center}
\caption{The action of the map $\gamma_0$: Vortices are moved on equidistant trajectories around the hole $\Omega_1$.}
\end{figure}

\begin{proof2}
Let $\gamma\in\Gamma_0$, and let $\tilde H:\IT^N\times[0,1]\to\CF_N\Omega$ be a homotopy connecting $\gamma$ and $\gamma_0$. Setting \[
r:\Omega\ni z\mapsto \frac{z}{|z|}\in S^1,
\]\[
\Psi:=r\times\dots\times r:\CF_N\Omega\to\IT^N,
\]
$\Psi$ is well--defined and continuous, and the assertion is equivalent to $\bar e\in\Psi(\gamma(\IT^N))$, where
$\bar e=\left(e^{\frac{2\pi ij}{N}}\right)_{j\in\{1,\dots,N\}}\in\IT^N$.
Now for every $t\in[0,1]$, the map $f_t:=\Psi\circ \tilde H(\cdot,t)$ induces a homomorphism \[
f_*:\IZ\cong H_N(\IT^N)\to H_N(\IT^N)\cong\IZ,
\]
in singular homology which is independent of $t$. Since $h_j$ is a homeomorphism onto its image and $h_i\simeq h_j$ for $i,j\in\{1,\dots,N\}$, the map $r\circ h_j:S^1\to S^1$ has winding number
$\pm 1$, hence induces an isomorphism $(r\circ h_j)_*=(r\circ h_1)_*:h_1(S^1)\to h_1(S^1)$. Now if $\{S^1\}\in H_1(S^1)\cong\IZ$ is a generator, $\{\IT^N\}=\{S^1\}\times\dots\times\{S^1\}$ is a generator of $H_N(\IT^N)$ and we compute\[
f_*(\{\IT^N\})=(\Psi\circ\gamma_0)_*(\{S^1\}\times\dots\times\{S^1\})
=\left((r\circ h_1)\times\dots\times(r\circ h_N)\right)_*(\{S^1\}\times\dots\times\{S^1\})\]\[
=(r\circ h_1)_*(\{S^1\})\times\dots\times(r\circ h_1)_*(\{S^1\})
=\pm\{S^1\}\times\dots\times\{S^1\}=\pm\{\IT^N\},
\]
hence $f_*$ is an isomorphism. Now if $\bar e\not\in\Psi(\gamma(\IT^N))$, the isomorphism $f_*$ factorizes over $H_N(\IT^N\setminus\{\bar e\})$, that is we have a commutative diagram
\begin{diagram}
 \IZ\cong H_N(\IT^N) & & \rTo_{\cong}^{f_*} & & H_N(\IT^N)\cong\IZ\\
  & \rdTo  & & \ruTo_{j_*} & \\
& & H_N(\IT^N\setminus\{\bar e\}) & &
\end{diagram}
where $j:\IT^N\setminus\{\bar e\}\into\IT^N$. Further we have the exact sequence 
\begin{diagram}
H_N(\IT^N\setminus\{\bar e\}) &\rTo^{j_*} & H_N(\IT^N) & \rTo^\cong & H_N(\IT^N,\IT^N\setminus\{\bar e\}).
\end{diagram}
The restriction homomorphism to the right is an isomorphism since $\IT^N$ is a compact orientable connected $N$--dimensional manifold, see for example \cite{tomdieck}, chapter V, theorem 12.1.
Since the sequence is exact, we conclude that the homomorphism $j_*$ is trivial, which is a contradiction and the proof of \ref{intersection_pi1} is complete.
\end{proof2}

Equipped with the results of lemmata \ref{starbound} and \ref{intersection_pi1} we can now conclude via \ref{deformation} that there is a point $z\in\CF_N\Omega$, such that
\[\lim\limits_{t\to t^+(z)}H_\Omega(\phi(t,z))<\infty.\]
Hence we are in the position to apply proposition \ref{nocollisions}, which implies that $H_\Omega$ has to have a critical point, and therefore the proof of theorem \ref{notsimplyconnected} is finished as well.

\bibliographystyle{amsplain}
\enlargethispage{1cm}
\bibliography{literatur}

\providecommand{\bysame}{\leavevmode\hbox to3em{\hrulefill}\thinspace}
\providecommand{\MR}{\relax\ifhmode\unskip\space\fi MR }
\providecommand{\MRhref}[2]{%
  \href{http://www.ams.org/mathscinet-getitem?mr=#1}{#2}
}
\providecommand{\href}[2]{#2}
\begin{thebibliography}{10}

\bibitem{arefetal}
H.~Aref, P.K. Newton, M.~A. Stremler, T.~Tokieda, and D.~Vainchtein,
  \emph{{Vortex crystals}}, {Adv. Appl. Mech.} \textbf{{39}} ({2003}), {1--79}.

\bibitem{barry-etal}
Anna~M. Barry, Glen~R. Hall, and C.Eugene Wayne, \emph{Relative equilibria of
  the (1+n)-vortex problem}, Journal of Nonlinear Science \textbf{22} (2012),
  no.~1, 63--83.

\bibitem{bartolucci-pistoia:2007}
D.~Bartolucci and A.~Pistoia, \emph{{Existence and qualitative properties of
  concentrating solutions for the sinh-Poisson equation}}, {IMA J. Appl. Math.}
  \textbf{{72}} ({2007}), {706--729}.

\bibitem{bartsch-dai}
T.~Bartsch and Q.~Dai, \emph{{Periodic solutions of the $N$-vortex Hamiltonian
  system in planar domains}}, {Preprint, arXiv:1403.4533v2}, 2014.

\bibitem{ba-personal}
T.~Bartsch and A.~Pistoia, \emph{{Critical points of the N-vortex Hamiltonian
  in bounded planar domains and steady state solutions of the incompressible
  Euler equations}}, {SIAM J. Math. Anal.} (to appear).

\bibitem{ba-pi-we}
T.~Bartsch, A.~Pistoia, and T.~Weth, \emph{{N-vortex equilibria for ideal
  fluids in bounded planar domains and new nodal solutions of the sinh-Poisson
  and the Lane-Emden-Fowler equations}}, {Comm. Math. Phys.} \textbf{{297}}
  ({2010}), {653--686}.

\bibitem{ba-pi-we-err}
\bysame, \emph{{Erratum: N-vortex equilibria for ideal fluids in bounded planar
  domains and new nodal solutions of the sinh-Poisson and the Lane-Emden-Fowler
  equations}}, {Comm. Math. Phys.} \textbf{{33}} ({2015}), {1107}.

\bibitem{cao-liu-wei-1}
D.~Cao, Z.~Liu, and J.~Wei, \emph{Regularization of point vortices pairs for
  the euler equation in dimension two}, Arch. Rational Mech. Anal. (2013),
  1--39.

\bibitem{cao-liu-wei-2}
\bysame, \emph{Regularization of point vortices pairs for the euler equation in
  dimension two}, Arch. Rational Mech. Anal. \textbf{212} (2014), 179--217.

\bibitem{dai-diss}
Q.~Dai, \emph{{Periodic solutions of the N point-vortex problem in planar
  domains}}, Ph.D. thesis, Justus-Liebig-Universit\"at, Otto-Behaghel-Str. 8,
  35394 Gießen, 2014.

\bibitem{delpino-etal:2005}
M.~del Pino, M.~Kowalczyk, and M.~Musso, \emph{{Singular limits in
  Liouville-type equations}}, {Calc. Var. Part. Diff. Equ.} \textbf{{24}}
  ({2005}), {47--81}.

\bibitem{esposito-etal:2005}
P.~Esposito, M.~Grossi, and A.~Pistoia, \emph{{On the existence of blowing-up
  solutions for a mean field equation}}, {Ann. Inst. H. Poincar\'e, Anal. non
  lin.} \textbf{{22}} ({2005}), {227--257}.

\bibitem{flucher}
M.~Flucher, \emph{{Variational Methods with Concentration}}, {Birkhäuser},
  {Basel}, {1999}.

\bibitem{flucher-gustafsson}
M.~Flucher and B.~Gustafsson, \emph{{Vortex motion in two--dimensional
  hydrodynamics}}, TRITA-MAT-1997-MA-02 ({1997}).

\bibitem{gilbarg-trudinger}
D.~Gilbarg and N.~S. Trudinger, \emph{Elliptic partial differential equations
  of second order}, Springer, 1997.

\bibitem{gustafsson}
B.~Gustafsson, \emph{{On the convexity of a solution of Liouville's equation}},
  {Duke Math.J.} \textbf{{60}} ({1990}), {303--311}.

\bibitem{diss}
C.~Kuhl, \emph{{Stationary solutions to the N-vortex problem}}, Ph.D. thesis,
  Justus-Liebig-Universit\"at, Otto-Behaghel-Str. 8, 35394 Gießen, 2013.

\bibitem{sym_part}
\bysame, \emph{{Symmetric equilibria for the N-vortex problem}}, {Journal of
  Fixed Point Theory and Applications} (to appear).

\bibitem{mallier-maslowe}
R.~{Mallier} and S.~A. {Maslowe}, \emph{{A row of counter-rotating vortices}},
  Physics of Fluids \textbf{5} (1993), 1074.

\bibitem{marchioro-pulvirenti}
C.~Marchioro and M.~Pulvirenti, \emph{{Mathematical Theory of Incompressible
  Nonviscous Fluids}}, {Springer}, {New York}, {1994}.

\bibitem{newton}
P.~K. Newton, \emph{{The N--vortex problem}}, {Springer}, {Berlin}, {2001}.

\bibitem{newton2}
\bysame, \emph{$n$-vortex equilibrium theory}, Discr. Cont. Dyn. Syst.
  \textbf{19} (2007), 411--418.

\bibitem{roberts}
Gareth~E. Roberts, \emph{Stability of relative equilibria in the planar
  $n$-vortex problem}, SIAM Journal on Applied Dynamical Systems \textbf{12}
  (2013), no.~2, 1114--1134.

\bibitem{saffman}
P.~G. Saffman, \emph{{Vortex Dynamics}}, {Cambridge University Press},
  {Cambridge}, {1992}.

\bibitem{schechter}
M.~Schechter, \emph{{Linking Methods in Critical Point Theory}}, {Birkhäuser},
  {1999}.

\bibitem{tomdieck}
T.~{tom Dieck}, \emph{{Topologie}}, {2} ed., {de Gruyter}, {2000}.

\end{thebibliography}

\end{document}